\documentclass[11pt,a4paper]{article}
\usepackage[margin=27mm]{geometry}
\usepackage{amsmath,amssymb,amsthm}
\IfFileExists{lmodern.sty}{%
  \usepackage[T1]{fontenc}
  \usepackage{lmodern}
}{}
\usepackage{color}
\definecolor{paperlink}{rgb}{0,0,0.45}
\usepackage[colorlinks=true,linkcolor=paperlink,citecolor=paperlink,urlcolor=paperlink]{hyperref}
\hypersetup{pdftitle={A Sharp Nodal Measure Bound for the Mixed Clamped--Navier Biharmonic Eigenvalue Problem in Dimensions n >= 2},pdfsubject={Local all-centers clamped doubling and optimal quarter-power nodal measure},pdfkeywords={biharmonic, nodal measure, smooth boundary, clamped, Navier, sharp estimate}}
\numberwithin{equation}{section}
\newtheorem{theorem}{Theorem}[section]
\newtheorem{lemma}[theorem]{Lemma}
\newtheorem{proposition}[theorem]{Proposition}
\newtheorem{corollary}[theorem]{Corollary}

\theoremstyle{remark}
\newtheorem{remark}[theorem]{Remark}
\newcommand{\R}{\mathbb R}
\newcommand{\plateComplex}{\mathbb C}
\newcommand{\HH}{\mathcal H}
\newcommand{\Z}{\mathcal Z}
\newcommand{\dd}{\,\mathrm d}
\newcommand{\dist}{\operatorname{dist}}

\newcommand{\diver}{\operatorname{div}}
\newcommand{\norm}[2][]{\left\|#2\right\|_{#1}}

\newcommand{\GN}{\Gamma_{\mathrm N}}
\newcommand{\GC}{\Gamma_{\mathrm C}}
\title{A Sharp Nodal Measure Bound for the Mixed Clamped--Navier Biharmonic Eigenvalue Problem in Dimensions $n\ge 2$}
\author{\begin{tabular}{@{}c@{\hspace{3em}}c@{}}
    Hairong Liu & Long Tian\\[-1pt]
    \small\href{mailto:hrliu@njust.edu.cn}{\nolinkurl{hrliu@njust.edu.cn}} &
    \small\href{mailto:tianlong19850812@163.com}{\nolinkurl{tianlong19850812@163.com}}
  \end{tabular}\\[0.25em]
  \small School of Mathematics and Statistics, Nanjing University of Science and Technology,\\[-0.1em]
  \small Nanjing 210094, China.}
\date{Revised November 18, 2026}
\begin{document}
\maketitle
\begin{abstract}

We estimate the interior nodal set of real eigenfunctions of $\Delta^2u=\lambda u$ on a bounded connected Euclidean domain. Each boundary component is either of class $C^\infty$ with clamped conditions or of class $C^{2,\alpha}$ and subject to Navier conditions. We obtain $\mathcal H^{n-1}(\{u=0\}\cap\Omega)\le C_\Omega\lambda^{1/4}$ for $n\ge2$. A spectral inequality on auxiliary clamped collars supplies fixed-scale interior mass. On clamped patches we localize the Carleman estimates of Liu, Tian and Yang. On Navier patches a single harmonic function with trace $u$ controls the two boundary data required for subdivision at every spatial scale. Summing an index linear in local doubling and frequency gives the stated bound. An annulus with one component of each type realizes its order of growth.

\end{abstract}

\section{Introduction and precise statements}\label{m:intro}
Let $\Omega\subset\R^n$, $n\ge2$, be bounded and connected. Its boundary is a finite disjoint union of compact embedded hypersurfaces. Write
\[
 \partial\Omega=\GN\mathbin{\dot\cup}\GC,
\]
where each part is a union of complete components. Our main assumptions are
\begin{equation}\label{m:smooth-geometry}
 \GC\in C^\infty,\qquad \GN\in C^{2,\alpha},\qquad 0<\alpha<1.
\end{equation}
All chart radii and coordinate bounds are fixed. If both parts are nonempty, their compactness and separation give $\dist(\GC,\GN)>0$. There are no interfaces where the boundary condition switches along one component, and no corners. No convexity is assumed.

We study real nonzero eigenfunctions of
\begin{equation}\label{m:problem}
\Delta^2u=\lambda u\quad\text{in }\Omega,
\qquad u=\Delta u=0\text{ on }\GN,
\qquad u=\partial_\nu u=0\text{ on }\GC.
\end{equation}
The variational meaning is
\begin{equation}\label{m:weak}
\int_\Omega\Delta u\Delta\varphi\dd x
=\lambda\int_\Omega u\varphi\dd x\quad(\varphi\in\mathcal V),
\qquad
\mathcal V=\{w\in H^2(\Omega)\cap H^1_0(\Omega):\partial_\nu w|_{\GC}=0\}.
\end{equation}

The normal trace is continuous on $H^2(\Omega)$, so $\mathcal V$ is closed in $H^2\cap H^1_0$. On this Dirichlet space the global $H^2$ estimate and Poincar\'e's inequality give $\|w\|_{H^2}\le C_\Omega\|\Delta w\|_2$. Since $C_c^\infty(\Omega)\subset\mathcal V$ and $\mathcal V\hookrightarrow L^2(\Omega)$ is compact, the closed form defines a positive self-adjoint operator with compact resolvent. More precisely, if $\lambda_1^{\rm D}$ is the first Dirichlet Laplace eigenvalue, then
\[
 \lambda_1^{\rm D}\|w\|_2^2\le\|\nabla w\|_2^2
 =-\langle w,\Delta w\rangle\le\|w\|_2\|\Delta w\|_2.
\]
Thus every plate eigenvalue satisfies $\lambda^{1/2}\ge\lambda_1^{\rm D}>0$. The Navier trace in \eqref{m:problem} follows from Lemma~\ref{lem:transposition}. Set

\[
 k=\lambda^{1/4},\qquad K=1+k,
 \qquad\Z(u)=\{x\in\Omega:u(x)=0\}.
\]
Only the interior nodal set is counted, without multiplicity. The estimates are invariant under multiplication of $u$ by a nonzero constant; when needed we normalize $\norm[L^\infty(\Omega)]{u}=1$.

\begin{theorem}[Sharp nodal bound with smooth clamped components]\label{m:sharp}
Under \eqref{m:smooth-geometry} and the separated-component assumptions above, every real nonzero variational eigenfunction of \eqref{m:problem} satisfies
\begin{equation}\label{m:sharp-bound}
 \HH^{n-1}(\Z(u))\le C_\Omega\lambda^{1/4}.
\end{equation}
The constant depends on the dimension, the fixed boundary atlases, their separation and the geometry of $\Omega$. It is uniform in $u$ and $\lambda$. The exponent is attained in order along a sequence on a fixed mixed annulus, as shown in Proposition~\ref{m:mixed-sharpness}.
\end{theorem}

\begin{corollary}[Finitely many separated holes]\label{m:holes-theorem}
Let $D,H_1,\ldots,H_m$ be bounded domains, with $\overline H_j\subset D$ and pairwise disjoint $\overline H_j$. Suppose $\Omega=D\setminus\bigcup_{j=1}^m\overline H_j$ is connected. Give each of its finitely many boundary components either boundary condition in \eqref{m:problem}, with the corresponding regularity in \eqref{m:smooth-geometry}. Then \eqref{m:sharp-bound} holds.
\end{corollary}
\begin{proof}
The compact disjoint components have positive minimum separation. Their one-sided collars admit a common fixed width, and connectedness supplies finite chains of interior balls between the finitely many observation patches. This gives the geometry of Theorem~\ref{m:sharp}. Its constant may depend on the number and separation of the holes, the charts and the widths of the connecting passages. The assertion concerns a fixed domain.
\end{proof}

\paragraph{Relation to earlier work.}
In the pure Navier case the plate realization is $A^2$, where $A=-\Delta_{\rm D}>0$. The spectral theorem gives $Au=\sqrt\lambda\,u$, so the Dirichlet nodal bound of Logunov, Malinnikova, Nadirashvili and Nazarov \cite{LMNN} gives \eqref{m:sharp-bound}. When a clamped component is present, $\Delta u$ need not vanish on it, and this reduction ceases to apply. The interior analytic counting argument follows the approach of Donnelly--Fefferman and Lin \cite{DF,Lin}. Lin--Zhu \cite{LZ} obtained the quarter-power bound for analytic clamped domains. Tian--Yang \cite{TY} treated specified piecewise analytic geometries. Their spectral convention $\Delta^2u=\Lambda^2u$ makes their $\sqrt\Lambda$ equal to $\lambda^{1/4}$ here. For the growth theory of fourth-order equations, see Zhu's doubling and nodal estimates \cite{Zhu19} and his subsequent polynomial measure bounds \cite{Zhu26}. Morassi, Rosset and Vessella \cite{MRV} give boundary doubling for the planar isotropic Kirchhoff--Love equation with clamped data.

Liu, Tian and Yang \cite{LTY} treat the fully clamped smooth Euclidean problem. The annular and half-space Carleman estimates of \cite{LTY} apply to clamped test functions. Section~\ref{m:local-section} uses them on patches meeting a single clamped component and derives the doubling and subdivision estimates for the mixed eigenfunction there. The mixed problem adds the transfer of fixed-scale mass between different boundary components and a Navier estimate at arbitrary local frequency. The latter uses the same harmonic lift on every child of a boundary cube. In particular its large-index threshold has size $1+k\ell$ on a cube of side $\ell$.

\paragraph{Proof of the mixed estimate.}
Section~\ref{m:regularity} derives the Navier trace from the form, obtains lower mass bounds on fixed interior balls, and states the interior nodal estimate. Section~\ref{m:local-section} obtains the local growth bounds and the clamped child estimates. Appendices~\ref{sec:tools}--\ref{sec:recursion} give the Navier reflection, the two-trace estimates, and the controlled children. Section~\ref{m:sharp-section} sums the two boundary contributions and constructs the sharp sequence on an annulus.

\section{Fixed-scale observation and interior estimates}\label{m:regularity}

\begin{lemma}[Recovery of the Navier trace]\label{lem:transposition}
Let $P'\Subset P$ be ambient neighborhoods of a pure Navier patch, with $\overline P\cap\Gamma_{\rm C}=\varnothing$. Put $D_P=P\cap\Omega$. Suppose $v,f\in L^2(D_P)$ and $\int_{D_P}v\Delta\psi=\int_{D_P}f\psi$ for $H^2$ tests of zero Dirichlet trace supported away from the artificial boundary. Then $v|_{D_{P'}}\in H^2(D_{P'})$, its trace on $P'\cap\partial\Omega$ is zero, and
\[
 \|v\|_{H^2(D_{P'})}\le C_{P',P,\Omega}
       (\|v\|_{L^2(D_P)}+\|f\|_{L^2(D_P)}).
\]
In particular the assertion applies to $v=\Delta u$ and $f=\lambda u$ in \eqref{m:weak}.
\end{lemma}
\begin{proof}
Choose smooth real cutoffs $\eta,\zeta$ compactly supported in $P$, with $\zeta=1$ on $P'$ and $\eta=1$ near $\operatorname{supp}\zeta$. For $\phi\in H^2(\Omega)\cap H^1_0(\Omega)$ insert $\eta\phi$ in the identity. Extending the cutoff quantities by zero on the rest of $\Omega$ gives
\[
 \int_\Omega (\eta v)\Delta\phi
 =\int_\Omega (\eta f-v\Delta\eta)\phi-2\int_\Omega v\nabla\eta\cdot\nabla\phi .
\]
The right side defines $G\in H^{-1}(\Omega)$ with norm bounded by the stated data. Let $w\in H^1_0(\Omega)$ solve $\Delta w=G$. Then $\int(\eta v-w)\Delta\phi=0$ for all such $\phi$. The Dirichlet Laplacian maps $H^2\cap H^1_0$ onto $L^2$ on the $C^{1,1}$ domain $\Omega$, by the Dirichlet $H^2$ estimate \cite[Theorem~9.15, $p=2$]{GT}. Hence $w=\eta v$ a.e. In particular the localized trace is zero and $\|\eta v\|_{H^1}\le C\|G\|_{H^{-1}}$. Now $\zeta v\in H^1_0(\Omega)$ and
\[
 \Delta(\zeta v)=\zeta f+2\nabla\zeta\cdot\nabla(\eta v)+(\Delta\zeta)\eta v\in L^2(\Omega).
\]
The same Dirichlet estimate gives the asserted $H^2$ bound. For a plate eigenfunction the tests $\eta\phi$ belong to $\mathcal V$, which proves the last assertion.
\end{proof}

\begin{proposition}[A cutoff consequence of clamped spectral observation]\label{m:cutoff}
Let $D$ be a fixed smooth one-sided collar of one component of $\GC$, bounded on its other side by a smooth artificial interior hypersurface. There are a cutoff $\chi$ equal to one near the physical component and a fixed set $A\Subset\Omega$ such that every solution of \eqref{m:problem} satisfies
\begin{equation}\label{m:cutoff-bound}
\norm[L^2(D)]{\chi u}\le Ce^{CK}\norm[L^2(A)]{u}.
\end{equation}
No analyticity of $D$ is assumed.
\end{proposition}
\begin{proof}
Choose $\chi$ to vanish in a neighborhood of the artificial inner boundary. Then $v=\chi u$ belongs to the operator domain of the fully clamped realization $P_D$ of $\Delta^2$ on $D$, and
\[
(P_D-\lambda)v=f,\qquad f=[\Delta^2,\chi]u.
\]
Local regularity near the smooth clamped component and interior regularity justify these identities in $L^2$. The support of $f$ is a fixed interior strip; choose a slightly larger $A\Subset\Omega$ containing it.

Let $\Pi$ be the spectral projection of $P_D$ onto eigenvalues at most $2(1+\lambda)$. On the complementary space the distance to $\lambda$ is at least $1+\lambda$, so
\[
\norm[L^2(D)]{(I-\Pi)v}\le(1+\lambda)^{-1}\norm[L^2(D)]{f}.
\]
Choose a fixed nonempty open $\omega\subset D$ on which $v=0$. The clamped spectral inequality of Le Rousseau--Robbiano \cite[Theorem 1.3]{LR} gives
\[
\norm[L^2(D)]{\Pi v}\le Ce^{CK}\norm[L^2(\omega)]{\Pi v}
=Ce^{CK}\norm[L^2(\omega)]{(I-\Pi)v}.
\]
Hence $\norm[L^2(D)]{v}\le Ce^{CK}(1+\lambda)^{-1}\norm[L^2(D)]{f}$. The commutator has order at most three, and buffered interior elliptic estimates give $\norm[L^2(D)]{f}\le C(1+\lambda)\norm[L^2(A)]{u}$. This proves \eqref{m:cutoff-bound}. Each component can be treated separately if the chosen collection of collars is disconnected.
\end{proof}

\begin{lemma}[Finite-order bounds and energy]\label{m:finite}
There are an integer $M$ and a constant $C$ such that
\begin{equation}\label{m:finite-bound}
\norm[L^\infty(\Omega)]{u}+\norm[L^\infty(\Omega)]{\Delta u}
\le C(1+\lambda)^M\norm[L^2(\Omega)]{u}.
\end{equation}
Both fields are continuous on $\overline\Omega$, and
\begin{equation}\label{m:energy}
\norm[L^2(\Omega)]{\Delta u}=k^2\norm[L^2(\Omega)]{u}.
\end{equation}
Consequently the normalization implies $\norm[L^2(\Omega)]{u}\ge e^{-CK}$.
\end{lemma}
\begin{proof}
Testing \eqref{m:weak} with $u$ gives \eqref{m:energy}. On a Navier patch let $v=\Delta u$. The variational formulation yields $\Delta v=\lambda u$ with zero Dirichlet data there in the transposition sense; Lemma~\ref{lem:transposition} gives its $H^2$ regularity and trace. Apply local Dirichlet $W^{2,p}$ estimates to $\Delta u=v$ and $\Delta v=\lambda u$ on nested patches, starting with $p=2$. A finite number of Sobolev improvements gives an exponent $p>n$ and hence bounded continuous $u,v$. Each step introduces only a fixed power of $1+\lambda$. The required Dirichlet estimates hold on $C^{2,\alpha}$ patches. Clamped patches are smooth, and a finite number of the fourth-order boundary estimates \cite{ADN} gives the same conclusion for $u$ and its second derivatives. Interior patches are handled in the same way. A finite buffered cover gives \eqref{m:finite-bound}, using \eqref{m:energy} for the initial norm of $v$. Taking $\norm[L^\infty]{u}=1$ gives a polynomial lower bound for its $L^2$ norm; fixed powers of $1+\lambda$ are bounded by $e^{CK}$.
\end{proof}

We use the exact scalar harmonic lift
\begin{equation}\label{m:H}
H(x,t)=C_k(t)u(x)-S_k(t)\Delta u(x),\qquad
C_k(t)=\frac{\cosh kt+\cos kt}{2},\quad
S_k(t)=\frac{\cosh kt-\cos kt}{2k^2}.
\end{equation}
The positive spectral lower bound in Section~\ref{m:intro} keeps $k$ away from zero on a fixed domain. At $k=0$ the continuous extension is $S_0(t)=t^2/2$. Direct differentiation gives $\Delta_{x,t}H=0$ and $H(x,0)=u(x)$. On each fixed slab its supremum is at most $e^{CK}$, by Lemma~\ref{m:finite}. On $\GN\times\R$ it has zero Dirichlet trace. The full lifting and norm comparison statements appear in Appendix~\ref{sec:lift}.

\begin{lemma}[Uniform seeds on fixed interior balls]\label{m:seeds}
For every fixed ball $B\Subset\Omega$ there is $C_B>0$ such that
\begin{equation}\label{m:seed}
\norm[L^2(B)]{u}\ge e^{-C_BK}.
\end{equation}
For centers in a fixed compact subset of $\Omega$ and a fixed admissible radius, the constants can be chosen uniformly. The corresponding fixed balls centered at $(x,0)$ for the lift $H$ also have $L^2$ norm at least $e^{-CK}$.
\end{lemma}
\begin{proof}
Choose small fixed collars of the clamped components, inside the regions where the cutoffs in Proposition~\ref{m:cutoff} equal one. The complement is covered by finitely many fixed interior sets and smaller Navier boundary patches with larger buffers. From Lemma~\ref{m:finite}, one member of this finite cover carries $L^2$ norm of $u$ at least $e^{-CK}$.

If it is a clamped collar, Proposition~\ref{m:cutoff} implies an $e^{-CK}$ lower bound on a fixed interior observation strip $A$. A finite cover of $A$ gives a fixed interior ball with the same exponential type of lower bound. If it is already an interior set, a finite cover gives such a ball directly.

If it is a Navier patch, flatten and reflect $H$ oddly as in Lemma~\ref{lem:reflection}. It solves a scalar uniformly elliptic Lipschitz equation with coefficients independent of $k$. The trace inequality in the $t$ variable and elliptic Caccioppoli give
\[
\norm[L^2(\text{smaller spatial patch})]{u}
\le C\norm[L^2(\text{buffered reflected product patch})]{H}.
\]
The product patch on the right is chosen inside a larger reserved one. Its norm is therefore at least $e^{-CK}$, while the larger norm is at most $e^{CK}$. A fixed-ratio three-ball inequality, with an inner ball lying on the physical interior side and a middle ball containing this smaller patch, gives an $e^{-CK}$ lower bound on that inner ball. Equivalently a finite chain of fixed-radius reflected balls gives the same conclusion. All radii and chain lengths are fixed geometrically. Thus a fixed interior product ball carries lifted norm at least $e^{-CK}$.

There is no difficulty passing between scalar and lifted norms on fixed interior balls. The upper comparison \eqref{eq:large-parameter-upper}, which follows from interior Caccioppoli, bounds the lifted norm by $e^{CK}$ times the scalar norm on a fixed enlarged spatial ball. In the other direction, the harmonic trace estimate \eqref{eq:large-parameter-lower} bounds the scalar norm by a constant times a fixed enlarged lifted norm. Consequently all three cases provide a fixed interior lifted seed.

Connectedness of $\Omega$ permits a finite chain of interior product balls from each of these finitely many possible seeds to any given fixed target ball. Harmonic three-ball interpolation and the global upper bound $e^{CK}$ propagate the lower bound along the chain. Make its last lifted ball small enough that its enlarged spatial projection lies in the prescribed $B$, and apply the upper scalar-to-lift comparison once more. This proves \eqref{m:seed}. Uniformity on compact sets follows by a finite net, applying the argument to smaller balls contained in all admissible target balls. Finally a fixed lifted ball contains a smaller product cylinder, and the harmonic trace estimate on an even smaller spatial ball proves the asserted lower bound for $H$.
\end{proof}

The constants in \eqref{m:seed} depend on the chosen radius. Section~\ref{m:local-section} supplies the growth estimates as the radius tends to zero.

\subsection{Counting zeros in an interior ball}

\begin{lemma}[Interior analytic nodal estimate]\label{lem:analytic}
Let $w$ be a real, nonzero solution of $\Delta^2w=k^4w$ in $B_{L_1s}(x_0)\subset\R^n$. Choose $n+1$ anchor balls of radius $\epsilon_0s$ near the vertices of a fixed nondegenerate simplex in $B_s(x_0)$. If
\[
\norm[L^2(D_j)]{w}\geq e^{-D}\norm[L^2(B_{L_1s})]{w}
\quad(0\leq j\leq n),
\]
then
\begin{equation}\label{eq:analytic-measure}
\HH^{n-1}(\Z(w)\cap B_s(x_0))\leq Cs^{n-1}(1+ks+D).
\end{equation}
The constants $L_1,\epsilon_0,C$ depend only on the dimension and the fixed simplex geometry.
\end{lemma}
\begin{proof}
The complex-valued lift $\widetilde w(x,t)=\exp(e^{i\pi/4}kt)w(x)$ solves the parameter-independent strongly elliptic equation $(\Delta_x^2+\partial_t^4)\widetilde w=0$. Interior analytic estimates \cite{MN} give, for $L_1$ sufficiently large, a holomorphic extension of $w$ with
\[
\sup_{B_{4s}^{\plateComplex}(x_0)}|w^{\plateComplex}|
\leq Cs^{-n/2}e^{Cks}\norm[L^2(B_{L_1s})]{w}.
\]
Select $p_j\in D_j$ with $|w(p_j)|\geq cs^{-n/2}e^{-D}\norm[L^2(B_{L_1s})]{w}$. Restriction to a complex line through $p_j$ and Jensen's formula show that each ray from $p_j$ meets the nodal set in $B_s(x_0)$ at most $C(1+ks+D)$ times. The restriction is not identically zero because $w(p_j)\ne0$.

For radial projection from $p_j$ to the unit sphere, the tangential Jacobian on the top-dimensional strata of the zero set is
\[
J_j(x)=\frac{|\nu_Z(x)\cdot(x-p_j)|}{|x-p_j|^n}.
\]
The anchors have width at least $cs$ in every direction, and $|x-p_j|\leq Cs$. Therefore $\sum_jJ_j(x)\geq cs^{-(n-1)}$. The area formula \cite{Federer}, followed by the ray counts, proves \eqref{eq:analytic-measure}. Lower-dimensional strata have zero $\HH^{n-1}$ measure. The argument can first be applied to compact interior subsets, so it does not require a measure bound in advance.
\end{proof}

\section{Growth estimates near the two boundary types}\label{m:local-section}
\subsection{Local doubling on a clamped patch}
The clamped analysis uses two local inequalities from the fully clamped problem \cite[Theorem~3.1 and Lemma~4.3]{LTY}. We state them below with their weights, support conditions and parameter ranges. The first applies to a cutoff on a smooth auxiliary collar; the second applies to the flattened half-space. Their combination gives Theorem~\ref{m:local-doubling} for a solution defined only on the mixed-domain patch. Proposition~\ref{m:cutoff} supplies the fixed-scale observation used later to control its outer mass ratio.

Fix buffered neighborhoods $U_C^0\Subset U_C$ of a clamped boundary piece, with containment relative to a neighborhood of $\GC$, such that $U_C$ meets no Navier component. All artificial boundaries remain outside the supports used below. Write
\begin{equation}\label{m:J-local}
 J_x(r)=\int_{B_r(x)\cap\Omega}|u(y)|^2\dd y,
 \qquad D_L(x,r)=\log\frac{J_x(Lr)}{J_x(r)}.
\end{equation}
The denominators are positive: a nonzero solution is analytic in the connected interior, and cannot vanish on an interior open set. The same conclusion holds for balls centered on the physical boundary.

\begin{theorem}[Local all-centers clamped doubling]\label{m:local-doubling}
For every fixed sufficiently large $L$, there are $R_C>0$, $c_C\in(0,1)$ and $C_C>0$ such that
\begin{equation}\label{m:HC}
 D_L(x,r)\le C_C\bigl[1+D_L(x,R)+kR\bigr],
 \qquad 0<r\le c_CR,\quad 0<R\le R_C,
\end{equation}
whenever $x\in\overline\Omega\cap U_C^0$ and the fixed outer buffers, including $B_{LR}(x)$, lie in $U_C$. For the local statement let $u\in H^2$ solve $\Delta^2u=k^4u$, with $k\ge0$, on a connected open physical patch containing $\Omega\cap B_{LR}(x)$ and all cutoff neighborhoods, nonzero on that patch, with $u=\partial_\nu u=0$ on its physical face. Each sampled ball meets the interior in an open set, hence $J_x(s)>0$ at every radius in the formula by interior analyticity. All cutoffs vanish near the artificial boundary of the patch. The constants are uniform in the distance of the center to $\GC$, in $k$, and in the smaller radius. They depend on a finite number of fixed smooth chart bounds and the buffer geometry.
\end{theorem}

We first specify the two cited analytic inputs in local notation. For a physical dilation by $R$ about $x$, let $\Omega_{x,R}=R^{-1}(\Omega-x)$, $r=|y|$, and $\nu=kR$. The region $r<4$ is required to lie in the pure clamped buffer. Near the boundary let $d_R$ be inward signed distance in these scaled coordinates, and put
\begin{equation}\label{m:local-phase}
 T=\log(r_c/r),\qquad
 \Phi=T-a\log T+b\frac{d_R(y)-d_R(0)}{rT^2}.
\end{equation}
The fixed constants are chosen as in \cite[Section~3, (3.1)]{LTY}: first $a>0$, then a sufficiently large logarithmic offset in $r_c$, and finally $0<b\ll a$. For an interior pole with its required ball away from the boundary, omit the distance term. Shrinking $R_C$ makes the relevant tubular neighborhood uniform.

\begin{lemma}[Local annular analytic input]\label{m:annular-input}
There are a fixed finite exponent $M$ and a constant $C$ such that, for $0<\delta_*<1$,
\begin{equation}\label{m:annular-input-eq}
 \tau^3\norm[L^2]{e^{\tau\Phi}z}
 \le C T_*^M\norm[L^2]{r^4e^{\tau\Phi}(\Delta^2-\nu^4)z},
 \qquad \tau\ge C(1+\nu+T_*^M),\quad T_*=\log(r_c/\delta_*),
\end{equation}
for test functions smooth up to the physical face, with both clamped traces zero and support compactly contained, relative to $\overline\Omega_{x,R}$, in $\{\delta_*<r<3\}$. The enlarged annuli $r<4$ lie inside the pure patch. Here $x\in\overline\Omega$ and $\nu=kR$. The phase has a fixed gap
\begin{equation}\label{m:phase-gap}
 \inf_{0<r\le1}\Phi-\sup_{2<r<3}\Phi\ge g>0,
\end{equation}
and a bounded oscillation between radii with any prescribed fixed ratio. The constants are uniform in all centers, including interior centers tending to the physical boundary independently of $\delta_*$.
\end{lemma}
\begin{proof}[Identification and localization of the input]
Choose a smooth collar $D$ as in Proposition~\ref{m:cutoff}. Restrict $U_C^0$ and $R_C$ so that every boundary ball used here is separated from the artificial face of $D$ by more than its outer buffer radius. The physical parts of $D$ and $\Omega$ then agree on the support of $z$, and their inward signed distances agree throughout the required tubular neighborhood of $\Gamma_{\rm C}$. Extending $z$ by zero outside its support gives an admissible clamped test function on the rescaled collar. Theorem~3.1 of \cite{LTY}, applied on this fixed smooth domain with spectral parameter $\nu=kR$ and inner radius $\delta_*$, gives \eqref{m:annular-input-eq}. For an interior ball disjoint from the face one uses the radial phase of the same estimate. Each collar is fixed before choosing $R$ or $\delta_*$, so taking the minimum admissible radius and the maximum constant over the finite collection preserves uniformity as the pole approaches the boundary. The phase gap follows from $(T-a\log T)'\ge 3/4$ and $|d_R(y)-d_R(0)|\le r$ after reducing $b$. The same bounds control the oscillation on a fixed-ratio annulus.

For clarity, multiplication by any smooth cutoff preserves both physical traces:
\[
 (\chi u)|_{\GC}=0,\qquad
 \partial_\nu(\chi u)|_{\GC}
 =\chi\partial_\nu u+u\partial_\nu\chi=0.
\]
The radial cutoff and its derivatives vanish near the artificial boundary. Hence no new boundary term, and in particular no Navier boundary term, enters the local estimate. This is a localization of a test-function estimate; it does not extend $u$ as a global clamped eigenfunction on an auxiliary domain. Compactness of the finite pure clamped atlas supplies uniform constants.
\end{proof}

In the flat coordinates let $x=(0,d)$, $d\ge0$, $r=|y-x|$ and $\dd\gamma_x=r^{-n}\dd y$ on $\R^n_+$. Take $z$ smooth up to the plane, compactly supported in $\overline{\R^n_+}\setminus\{x\}$, with $z=\partial_nz=0$ on the plane. Its support may meet the plane. Denote the $L^2(\dd\gamma_x)$ norm by $\|\cdot\|_\gamma$ and put
\[
 \|w\|_{4,\sigma}^2=\sum_{j=0}^4\sigma^{2(4-j)}\|r^jD^jw\|_\gamma^2,\qquad
 A_i=\|r^{i-\sigma}D^iz\|_\gamma.
\]
Equation (4.6) of \cite[Lemma~4.3]{LTY} gives
$\sum_{i=0}^4\sigma^{1-i}A_i\le C\|r^{4-\sigma}\Delta^2z\|_\gamma$
for $\sigma\in\mathbb Z+\tfrac12$, $\sigma\ge4$, uniformly in $d$ and in the support radii. Since
$|D^b r^{-\sigma}|\le C_b\sigma^b r^{-\sigma-b}$ for $0\le b\le4$, Leibniz's rule gives
\begin{align*}
 \|r^jD^j(r^{-\sigma}z)\|_\gamma &\le C\sum_{i=0}^j\sigma^{j-i}A_i,\\
 \sigma^{-4}\|r^{-\sigma}z\|_{4,\sigma}
 &\le C\sum_{i=0}^4\sigma^{-i}A_i
 \le C\sigma^{-1}\|r^{4-\sigma}\Delta^2z\|_\gamma .
\end{align*}
In particular we may use the weaker bound
\begin{equation}\label{m:flat-input}
 \sigma^{-4}\|r^{-\sigma}z\|_{4,\sigma}
 \le C\|r^{4-\sigma}\Delta^2z\|_\gamma .
\end{equation}
The whole-space assertion of the same lemma applies to an interior patch. Keeping this weaker power permits the choice $\delta=\tau^{-5}$ below.

\begin{lemma}[Curved small-scale propagation]\label{m:small-local}
In a normalized flattened patch suppose
\[
 \mathcal L=\Delta_y^2+\sum_{|\beta|\le4}e_\beta(y)\partial^\beta,
 \qquad |e_\beta(y)|\le\mathfrak K r\ (|\beta|=4),\quad
 |e_\beta(y)|\le\mathfrak K\ (|\beta|\le3),
\]
where $r=|y-x|$, and suppose the coefficients and local clamped elliptic estimates have uniform fixed-order bounds. Fix $L_p\ge8$ and $0<a_0\le1/L_p$. Let $v$ solve $(\mathcal L-\nu^4)v=0$ on a connected buffered cap containing $B_{L_pa_0}(x)\cap\overline{\R^n_+}$, with clamped traces on its physical face and $J_c(t)>0$ at the radii used below. Set
\[
 J_c(t)=\int_{B_t(x)\cap\{y_n>0\}}|v|^2\dd y,
 \qquad N_0=\log\frac{J_c(L_pa_0)}{J_c(a_0)}.
\]
There are fixed $c_p,c_a,C_0>0$ such that, if
\begin{equation}\label{m:small-conditions}
 \sigma\in\mathbb Z+\tfrac12,\quad
 \sigma\ge\max\{4,C_0(1+N_0)\},\quad
 \mathfrak K a_0\sigma^4\le c_a,\quad \nu a_0\le c_a,
\end{equation}
then
\begin{equation}\label{m:small-conclusion}
 \log\frac{J_c(L_pt)}{J_c(t)}\le C_0\sigma,
 \qquad 0<t\le c_pa_0.
\end{equation}
The same result holds in an interior patch.
\end{lemma}
\begin{proof}
For a clamped test function supported in $0<r<a$, the perturbation terms obey
\begin{align*}
 \norm[L^2(\dd\gamma_x)]{r^{4-\sigma}(\mathcal L-\Delta^2)z}
 &\le C\mathfrak K a\norm[4,\sigma]{r^{-\sigma}z},\\
 \norm[L^2(\dd\gamma_x)]{\nu^4r^{4-\sigma}z}
 &\le \nu^4a^4\sigma^{-4}\norm[4,\sigma]{r^{-\sigma}z}.
\end{align*}
The derivatives of $r^{-\sigma}$ introduce only the powers of $\sigma$ already present in this norm. Absorbing these terms in \eqref{m:flat-input} gives
\begin{equation}\label{m:curved-singular}
 \sigma^{-4}\norm[4,\sigma]{r^{-\sigma}z}
 \le C\norm[L^2(\dd\gamma_x)]{r^{4-\sigma}(\mathcal L-\nu^4)z}
\end{equation}
when $\mathfrak K a\sigma^4$ and $\nu^4a^4$ are sufficiently small. Decrease $c_a$ to allow the support radius $a=3a_0$.

Choose an outer cutoff equal to one on $B_{2a_0}$ and zero outside $B_{3a_0}$, and an inner cutoff zero on $B_{t/4}$ and one outside $B_{t/2}$. Their product $\chi$ preserves the clamped traces. Applying \eqref{m:curved-singular} to $\chi v$, and using scaled local elliptic estimates on the cutoff annuli, gives
\begin{equation}\label{m:small-cutoff}
 \norm[L^2(\dd\gamma_x)]{r^{-\sigma}\chi v}
 \le C\left(4^\sigma t^{-\sigma-n/2}J_c(t)^{1/2}
 +2^{-\sigma}a_0^{-\sigma-n/2}J_c(L_pa_0)^{1/2}\right).
\end{equation}
These elliptic constants are independent of $\sigma,t,a_0$ because $\nu a_0$ is bounded. Every commutator term has at most three derivatives and the powers of the annular radius cancel after multiplication by $r^4$.

If $J_c(t)\ge J_c(a_0)/2$, monotonicity gives
$J_c(L_pt)/J_c(t)\le2e^{N_0}$ for $t\le a_0$.
Otherwise $\chi=1$ on $B_{a_0}\setminus B_t$, and the left side of \eqref{m:small-cutoff} is at least
$2^{-1/2}a_0^{-\sigma-n/2}J_c(a_0)^{1/2}$.
The outer term is at most $C2^{-\sigma}e^{N_0/2}$ times this quantity and is absorbed by the first condition in \eqref{m:small-conditions}. Retaining $B_{L_pt}\setminus B_t$ on the left, for $c_p\le1/L_p$, now gives
\[
 \frac{J_c(L_pt)}{J_c(t)}\le1+CL_p^n(4L_p)^{2\sigma}.
\]
This proves \eqref{m:small-conclusion} after increasing $C_0$. No extra annular observability input has been used.
\end{proof}

\begin{proof}[Proof of Theorem~\ref{m:local-doubling}]
All geometric constants are fixed before the final large numerical constant below. At a center within distance $4R$ of the physical boundary, use Fermi coordinates. A constant tangential linear change makes the tangential metric the identity at the pole and preserves the flattened boundary plane. The normal metric coefficient is one and the mixed normal--tangential coefficients vanish. Thus, after scaling by $R$, the pulled-back bi-Laplacian has the form in Lemma~\ref{m:small-local}, with
\begin{equation}\label{m:curvature-K}
 \mathfrak K\le C_{\rm geom}R.
\end{equation}
Normalization at the pole, rather than only at its projection on the boundary, ensures an $O(r)$ principal error even for a pole above the boundary. At an interior center use Euclidean coordinates and $\mathfrak K=0$.

Let $G\ge1$ be a uniform bi-Lipschitz constant, and let $m_-,m_+>0$ bound the volume density. Define the physical scaled mass
\[
 J(t)=R^{-n}J_x(Rt),\qquad N=D_L(x,R),\quad \nu=kR.
\]
For the coordinate mass $J_c$ we have
\begin{equation}\label{m:mass-compare}
 m_-J(t/G)\le J_c(t)\le m_+J(Gt).
\end{equation}
Set $L_p=G^2L$ and $\Gamma=G^2L_p=G^4L$.

Apply Lemma~\ref{m:annular-input} with an inner cutoff vanishing below $\rho/4$ and equal to one above $\rho/2$, and an outer cutoff equal to one below $2$ and zero from $11/4$ onwards. Local elliptic estimates give
\begin{equation}\label{m:annular-cutoff}
 \norm[L^2]{e^{\tau\Phi}\chi u_R}
 \le\Pi\left(e^{\tau\phi_i^+(\rho)}J(\rho)^{1/2}
       +e^{\tau\phi_o^+}J(4)^{1/2}\right),
 \qquad 1\le\Pi\le CT_*^p(1+\tau+\nu)^q.
\end{equation}
Here $p,q$ are fixed finite exponents and $\phi_i^+$ and $\phi_o^+$ are the phase suprema on the inner and outer cutoff annuli. There is no negative power of $\rho$ in $\Pi$: the terms $r^4[\Delta^2,\chi]u_R$ on the inner annulus are controlled by the norms of $\rho^jD^ju_R$, $0\le j\le3$, hence by a fixed polynomial in $1+\nu\rho$ times $J(\rho)^{1/2}$. The outer buffers lie inside $B_4$, and $J(4)\le J(L)=e^NJ(1)$.

Choose
\begin{equation}\label{m:tau-delta}
 \tau=A(1+N+\nu),\qquad \delta=\tau^{-5},\qquad
 \delta_* =\delta/16.
\end{equation}
Then $T_*=\log(16r_c)+5\log\tau$. For one fixed sufficiently large lower bound on $A$, the annular threshold holds and
\[
 \log\Pi\le\min\{1,g/4\}\tau,\qquad
 \sqrt2\Pi e^{N/2-g\tau}\le\tfrac12.
\]
These follow from $(\log\tau)^a/\tau\to0$ for each fixed finite $a$, and remain true if $A$ is increased. Require also $\delta\le\Gamma^{-1}$.

For $\delta\le\rho\le\Gamma^{-1}$, if $J(\rho)\ge J(1)/2$ then $J(\Gamma\rho)/J(\rho)\le2$. Otherwise the mass on $B_1\setminus B_\rho$ and \eqref{m:phase-gap} absorb the outer term in \eqref{m:annular-cutoff}. Retaining $B_{\Gamma\rho}\setminus B_\rho$ on the left and using the fixed-ratio phase oscillation gives
\begin{equation}\label{m:annular-range}
 \log\frac{J(\Gamma\rho)}{J(\rho)}\le C_1\tau,
 \qquad \delta\le\rho\le\Gamma^{-1}.
\end{equation}
The constant $C_1$ is independent of the final choice of $A$.

Now fix the constants $c_p,c_a,C_0$ from Lemma~\ref{m:small-local}, and choose
\[
 K_*\ge\max\{2G,2/(Gc_p)\},\qquad a_0=K_*\delta.
\]
Increase the lower bound on $A$ so that all coordinate balls through radius $L_pa_0$ are in the reserved chart, $a_0/G\le\Gamma^{-1}$, and $a_0\le1/L_p$. Since $a_0/G\ge\delta$, \eqref{m:mass-compare} and \eqref{m:annular-range} imply
\begin{align*}
 N_0:=\log\frac{J_c(L_pa_0)}{J_c(a_0)}
 &\le\log(m_+/m_-)+\log\frac{J(GL_pa_0)}{J(a_0/G)}\\
 &=\log(m_+/m_-)+\log\frac{J(\Gamma a_0/G)}{J(a_0/G)}
 \le C_2\tau.
\end{align*}
Choose the least half-integer satisfying $\sigma\ge\max\{4,C_0(1+N_0)\}$. Then
\begin{equation}\label{m:sigma-control}
 4\le\sigma\le C_\sigma\tau,
\end{equation}
where $C_\sigma$ is fixed before the final $A$. In view of \eqref{m:curvature-K},
\begin{equation}\label{m:matching}
 \mathfrak K a_0\sigma^4
 \le C_{\rm geom}R_CK_*C_\sigma^4\tau^{-1},
 \qquad \nu a_0\le K_*\tau^{-4}.
\end{equation}
One final fixed enlargement of $A$ makes both quantities at most $c_a$, as well as satisfying all earlier requirements. Thus Lemma~\ref{m:small-local} applies at every $0<t\le c_pa_0$.

Returning to physical balls using \eqref{m:mass-compare},
\[
 \frac{J(L\rho)}{J(\rho)}
 \le\frac{m_+}{m_-}\frac{J_c(GL\rho)}{J_c(\rho/G)}
 =\frac{m_+}{m_-}\frac{J_c(L_p\rho/G)}{J_c(\rho/G)}.
\]
This bounds the physical index by $C\tau$ for $0<\rho\le Gc_pa_0$, a range containing $(0,2\delta]$. The annular range covers $\delta\le\rho\le\Gamma^{-1}$, so the ranges overlap. Take $c_C=(2\Gamma)^{-1}$ and scale back by $R$ to obtain \eqref{m:HC}. Boundary centers follow either from the same estimates or by limits of interior centers: at each fixed positive radius the ball integrals vary continuously with the center, and their denominators stay positive. All steps have remained inside $U_C$.
\end{proof}

\begin{remark}[Order of constants and scope of the input]\label{m:local-dependence}
Fix the chart buffers, the ball-comparison constant $G$, and the ratio $L$ first. They determine $L_p,\Gamma$ and the constants in the annular and small-scale estimates. Next choose $K_*,C_\sigma$, and lastly $A$ in \eqref{m:tau-delta}. The exponent $5$ ensures that $a_0\sigma^4=O(\tau^{-1})$. Thus the logarithms from the annular threshold and the weight loss in the singular estimate are absorbed at the overlap of the two ranges.
\end{remark}

\subsection{Growth of the reflected Navier lift}\label{m:global}

\begin{lemma}[Spectral control of the lift near $\GN$]\label{m:navier-growth}
In a fixed finite atlas of reflected Navier patches, the maximal doubling index of the scalar harmonic lift is bounded by $CK$ at every smaller scale. Only the fixed interior seeds of Lemma~\ref{m:seeds} are needed for this conclusion.
\end{lemma}
\begin{proof}
Taking $\varphi=u$ in \eqref{m:weak} gives
\begin{equation}\label{m:navier-energy}
\norm[L^2(\Omega)]{\Delta u}=k^2\norm[L^2(\Omega)]{u}\le C_\Omega k^2.
\end{equation}
Define
\begin{equation}\label{m:global-lift}
H_k(x,t)=\frac{\cosh(kt)+\cos(kt)}2u(x)
-\frac{\cosh(kt)-\cos(kt)}{2k^2}\Delta u(x).
\end{equation}
The identities in Lemma~\ref{lem:lift} show that $H_k$ is harmonic and $H_k(x,0)=u(x)$. On a fixed slab its $L^2$ norm is at most $e^{CK}$ by \eqref{m:energy}; the formula has a continuous interpretation at $k=0$ as well. On $\GN\times\R$ it vanishes. The adapted flattening and odd reflection in Lemma~\ref{lem:reflection} yield a scalar uniformly elliptic divergence-form equation with Lipschitz coefficients independent of $k$.

Inside each smaller Navier chart choose a fixed physical ball $B_a(x_0)$ whose eightfold enlargement lies in $\Omega$. Lemma~\ref{m:seeds} supplies a point $x_1$ in that ball with $|u(x_1)|\ge e^{-CK}$. A fixed physical product ball around $(x_1,0)$ lies in the chart buffer. The harmonic mean-value estimate bounds this value by a constant times the $L^2$ norm on that product ball. Therefore a fixed inner lifted ball in the buffered chart has norm at least $e^{-CK}$, while its outer norm is at most $e^{CK}$.

Coordinate changes have bounded Jacobians, and a finite number of fixed-radius interior comparisons place this inner ball at the center used in the root estimate. The buffered frequency consequence Lemma~\ref{lem:frequency}(ii) now bounds all smaller lifted doubling indices, including centers in the reflected side, by $CK$. There are finitely many patches and their buffers stay away from $\GC$. No analytic continuation at the clamped boundary was used.
\end{proof}

\begin{corollary}[Global scalar doubling and logarithmic propagation]\label{m:global-scalar}
For the mixed eigenfunctions of Theorem~\ref{m:sharp}, there are fixed $r_D,C_D>0$ such that, for all $x\in\overline\Omega$ and $0<r\le r_D$,
\begin{equation}\label{m:global-L2-doubling}
 J_x(2r)\le e^{C_DK}J_x(r).
\end{equation}
After increasing $C_D$ if necessary, one also has, for the same centers and radii,
\begin{equation}\label{m:D}
\norm[L^\infty(B_{2r}(x)\cap\Omega)]{u}
\le e^{C_DK}\norm[L^\infty(B_r(x)\cap\Omega)]{u}.
\end{equation}
With the normalization $\|u\|_{L^\infty(\Omega)}=1$, there are $r_*\in(0,\min\{r_D,1\})$ and $C_P>0$ such that
\begin{equation}\label{m:Plog}
\norm[L^\infty(B_r(x))]{u}\ge
\exp\!\left[-C_PK\left(1+\log\frac{r_*}{r}\right)\right]
\quad\text{if }0<r\le r_*,\quad B_{8r}(x)\subset\Omega.
\end{equation}
The constants depend only on the fixed geometry and observation constants.
\end{corollary}
\begin{proof}
Near $\GC$, fix a small outer radius $R_*$ whose buffers remain in a pure clamped chart. The intersection $B_{R_*}(x)\cap\Omega$ contains an interior ball of a uniformly positive fixed radius with center in a fixed compact subset. By Lemma~\ref{m:seeds},
$J_x(R_*)\ge e^{-CK}$, while $J_x(LR_*)\le|\Omega|$.
Theorem~\ref{m:local-doubling} bounds all smaller fixed-ratio scalar indices by $CK$, including centers approaching the boundary.

Near $\GN$, Lemma~\ref{m:navier-growth} bounds the maximal doubling of the single reflected harmonic lift by $CK$. Its proof uses only Lemma~\ref{m:seeds}, not this corollary or the nodal-measure theorem. On comparable buffered balls, \eqref{eq:large-parameter-upper} and \eqref{eq:large-parameter-lower} compare lifted and scalar norms. Choose the inner lifted radius to be a fixed small multiple of $r$ so that its scalar buffer lies in $B_r(x)$, and choose the outer lifted radius to contain the trace neighborhood for $B_{2r}(x)$. A fixed number of lifted doublings then gives
\[
 \norm[L^2(B_{2r}(x)\cap\Omega)]{u}
 \le C e^{CK+Ckr}(1+(kr)^2)
          \norm[L^2(B_r(x)\cap\Omega)]{u}.
\]
The factors $r^{1/2}$ and $r^{-1/2}$ cancel, since the comparison radii have fixed ratios. For fixed $r\le r_D$, the remaining factors are absorbed into $e^{C'K}$. In the compact interior use the same lift, harmonic frequency estimates and the fixed seeds, with full balls. A finite cover proves \eqref{m:global-L2-doubling}.

The local finite-order estimates in Lemma~\ref{m:finite} and Appendix~\ref{sec:lift}, rescaled to $r$, give
\[
 \norm[L^\infty(B_{2r}(x)\cap\Omega)]{u}
 \le Cr^{-n/2}(1+(kr)^4)^M
       \norm[L^2(B_{4r}(x)\cap\Omega)]{u}
\]
for a fixed $M$. Apply \eqref{m:global-L2-doubling} twice and use
$\norm[L^2(B_r(x)\cap\Omega)]{u}\le Cr^{n/2}\norm[L^\infty(B_r(x)\cap\Omega)]{u}$.
This proves \eqref{m:D}. Fix $r_*<\min\{r_D,1\}/2$. Every $B_{r_*}(x)\cap\Omega$ contains an interior ball of uniformly positive radius with center in a fixed compact subset of $\Omega$, by the uniform interior ball condition. Lemma~\ref{m:seeds} and the volume bound give $\|u\|_{L^\infty(B_{r_*}(x)\cap\Omega)}\ge e^{-C K}$. For $r\le r_*$ choose the least integer $j\ge0$ such that $2^jr\ge r_*$. Iterating \eqref{m:D} at the radii $r,2r,\ldots,2^{j-1}r$ gives a lower bound of $e^{-(C+j C_D)K}$ on $B_r(x)\cap\Omega$. Now $j\le 1+\log_2(r_*/r)$, and $B_{8r}(x)\subset\Omega$ puts the entire ball $B_r(x)$ in the interior. This proves \eqref{m:Plog}.
\end{proof}

The lift \eqref{m:H} has different boundary traces on the two components. On $\Gamma_{\rm N}$ its zero trace permits odd reflection. On $\Gamma_{\rm C}$ it has trace $-S_k(t)\Delta u$, where $\Delta u$ is free. The clamped part of the preceding proof therefore uses Theorem~\ref{m:local-doubling}.

\subsection{Clamped boundary traces and child estimates}

\begin{lemma}[Spectral control of the clamped indices]\label{m:HC-seed}
Every sufficiently small fixed root cap near $\GC$ has a finite scalar maximal doubling index bounded by $CK$, uniformly over all its smaller scales and over centers in its fixed enlargement.
\end{lemma}
\begin{proof}
Fix a radius $R_*>0$ so small that all balls required by \eqref{m:HC} remain in the buffered clamped charts. Each $B_{R_*}(x)\cap\Omega$, with $x$ in a smaller clamped collar, contains an interior ball of a fixed smaller radius whose center stays in a fixed compact subset of $\Omega$. The uniform version of Lemma~\ref{m:seeds} gives
\[
J_x(R_*)\ge e^{-CK},\qquad J_x(LR_*)\le |\Omega|,
\]
after the normalization $\norm[L^\infty(\Omega)]{u}=1$. Consequently $D_L(x,R_*)\le CK$. Theorem~\ref{m:local-doubling} gives $D_L(x,r)\le CK$ for $r\le c_CR_*$. Choose all root sides below this scale. The choices are geometrical and independent of $k$.
\end{proof}

For a small-slope smooth graph use $F(y,t)=(y,g(y)+t)$. If $Q$ is a tangential cube of side $\ell$, put
\[
W(Q)=F(Q\times(0,\ell]),\qquad
E(Q)=F(16Q\times[0,16\ell]),
\]
enlarging the fixed numerical factor if needed. Reserve larger graph buffers. Set
\begin{equation}\label{m:clamped-index}
N_C(Q)=\sup\{D_L(x,r):x\in E(Q),\ 0<r\le\ell(Q)\},
\qquad \mu(Q)=k\ell(Q).
\end{equation}
For a tangential child $q\subset Q$ of ratio $A^{-1}$, the enlarged sets are nested, so
\begin{equation}\label{m:clamped-heredity}
N_C(q)\le N_C(Q),\qquad \mu(q)=\mu(Q)/A.
\end{equation}
All root scales are chosen so these sets and their PDE buffers lie in one clamped chart. The fixed large ratio $L$ is chosen to include the comparisons below before applying Theorem~\ref{m:local-doubling}.

\begin{lemma}[Clamped two-trace Cauchy propagation]\label{m:clamped-cauchy}
Normalize the base side to one and fix a central point $p$ of its graph. Let $W_{\rm out}=\Omega_{\rm loc}\cap B_{R_o}(p)$ in these rescaled variables, where $R_o$ is fixed large enough to cover the unit cap and the child comparisons. Suppose $\Delta^2v=\mu^4v$ there, with $v=\partial_\nu v=0$ on the physical face. Fix $a>0$ sufficiently small that the projection of $\Sigma=\partial\Omega_{\rm loc}\cap B_{4a}(p)$ lies in the central half of the base. Use the crossing ball $\mathbb B_{2a}(p,0)\subset\R^{n+1}$. Its physical projection lies in $B_{R_o}(p)$. Choose $B_*=B_\rho(y_0)$ and $I=(-h,h)$ so that $\overline{B_*\times I}$ lies strictly in its upper physical side inside $\mathbb B_a(p,0)$, and reserve a closed exterior ball in the lower side. These choices can be made uniformly on the small-slope graphs, with $y_0\in E(Q)$ and $0<\rho\le1$. There are $C>0$ and $\theta\in(0,1)$ such that
\begin{equation}\label{m:clamped-cauchy-eq}
\norm[L^2(B_*)]{v}\le Ce^{C\mu}M^{1-\theta}T^\theta,
\quad T=\norm[L^2(\Sigma)]{\Delta v}+\norm[L^2(\Sigma)]{\partial_\nu\Delta v},
\end{equation}
for every $M\ge\|v\|_{L^2(W_{\rm out})}$. The constants depend on the fixed chart bounds and the margins between these sets. All are chosen before subdivision, and $L$ in \eqref{m:clamped-index} is taken large enough that $W_{\rm out}\subset B_{L\rho}(y_0)$.
\end{lemma}
\begin{proof}
Put $\omega=e^{i\pi/4}$ and $V(x,t)=e^{\omega\mu t}v(x)$. The constant-coefficient operator $P=\Delta_x^2+\partial_t^4$ satisfies $PV=0$ and is strongly elliptic. Extend $V$ by zero across the physical cylindrical boundary. The clamped traces make this an $H^2$ extension locally. In a ball crossing a smaller central patch, integration by parts leaves only the two distributions determined by $\Delta v$ and $\partial_\nu\Delta v$ on the interface. With a fixed cutoff supported in the larger patch their sum $S$ has
\[
\norm[H^{-2}]{S}\le Ce^{C\mu}T.
\]
Here the dual space is that of $H_0^2$ on the crossing ball. Both the value and normal derivative of an $H^2$ test function have the required $L^2$ traces on the smooth interface, including the terms where the cutoff is differentiated.

The form $\int(|\Delta_x w|^2+|\partial_t^2w|^2)$ is coercive on $H_0^2$ of that ball, by its Fourier symbol and Poincar\'e's inequality. Solve $Pw=S$ there with $\norm[H^2]{w}\le Ce^{C\mu}T$. On a smaller crossing ball, the difference between the zero extension and $w$ solves $Ph=0$. Its norm on a fixed exterior ball is at most $Ce^{C\mu}T$, and its norm on the larger ball is at most $Ce^{C\mu}(M+T)$.

Here is the propagation argument for the homogeneous function $h$. Write $M_0$ for its $L^2$ norm on the larger ball and $\epsilon$ for its $L^2$ norm on the exterior observation ball. Cover a connected compact region joining a smaller observation ball to the target cylinder by finitely many overlapping balls of radius $r$, with consecutive centers at distance at most $r$. The compact region has positive distance from the boundary of the ball of homogeneity.
Choose $r$ small enough that the constant-coefficient interior analytic estimate \cite{MN} bounds the complex extension about every center on the complex ball of radius $4r$ by $C M_0$. If $\|h\|_{L^2(B_r(z))}\le b$, the same estimate on $B_r(z)$ bounds the extension on the concentric complex ball of radius $\eta r$ by $C_r b$, for a fixed $0<\eta<1$. All radii here are fixed independently of $\mu$. On every complex line through $z$, Hadamard's three-circles theorem at radii $\eta r,2r,4r$ gives
\[
 \sup_{B_{2r}(z)} |h| \le C_r M_0^{1-\gamma} b^\gamma,
 \qquad \gamma=\frac{\log 2}{\log(4/\eta)}.
\]
The next real ball of radius $r$ lies in $B_{2r}(z)$. Iteration along the finite cover bounds the $L^2$ norm on the target cylinder by $C M_0^{1-\theta}\epsilon^\theta$ for a fixed $\theta>0$. If $T\le M$, add back $w$ and use $|e^{\omega\mu t}|\ge e^{-C\mu}$ on this cylinder to obtain \eqref{m:clamped-cauchy-eq}; if $T>M$ the conclusion follows from the definition of $M$. The case $T=0$ follows by the same argument or by analytic unique continuation of the zero extension.
\end{proof}

\begin{lemma}[High-index clamped child]\label{m:clamped-high}
There are fixed $A_C\ge2$ and $M_C\ge1$ such that
\begin{equation}\label{m:clamped-drop}
N_C(Q)>M_C(1+\mu(Q))
\quad\Longrightarrow\quad
N_C(q)\le\tfrac12N_C(Q)
\end{equation}
for at least one of its $A_C^{n-1}$ boundary children.
\end{lemma}
\begin{proof}
Normalize $\ell(Q)=1$, put $N=N_C(Q)$ and $\mu=k\ell(Q)$, and let $M^2$ be the scalar mass on a fixed outer cap containing all comparisons. The fixed Cauchy geometry in Lemma~\ref{m:clamped-cauchy} gives an interior ball $B_\rho(y_0)$ with $y_0\in E(Q)$, $\rho\le1$. Choose $L$ initially so that $B_{L\rho}(y_0)$ contains the outer cap. Then
\begin{equation}\label{m:clamped-anchor}
\norm[L^2(B_\rho(y_0))]{u}\ge e^{-N/2}M.
\end{equation}
Suppose all children of a subdivision by $A$ have index greater than $N/2$. Choose witnesses $x_i\in E(q_i)$, $r_i\le A^{-1}$ with $D_L(x_i,r_i)>N/2$. Fix $B_0\ge c_C^{-1}$ large enough that the ball of radius $R_0=B_0/A$ about each witness contains buffered boundary regularity neighborhoods over that child.

Let $R_j=L^jR_0$, stopping before a fixed small $R_{\max}$ whose outer buffers remain in the parent chart and mass region. If $s$ is the number of steps, then $s\asymp\log A$ and $\sum_{j<s}R_j\le C$. For every step $r_i\le c_CR_j$, so \eqref{m:HC}, in these scaled coordinates, gives
\[
D_L(x_i,R_j)\ge\frac{N}{2C_C}-1-\mu R_j.
\]
Telescoping at this same center yields
\begin{equation}\label{m:clamped-small}
J_{x_i}(R_0)
\le C A^{p-cN}e^{C\mu}M^2.
\end{equation}
The constants in the powers are independent of $A$ once a fixed lower bound for $A$ is imposed. In particular the spectral terms sum geometrically and cost $C\mu$, not $C\mu\log A$.

Finite-order clamped boundary regularity on the buffered child neighborhoods gives
\[
\sup_{\Sigma_i}(|\Delta u|+|\partial_\nu\Delta u|)
\le C A^{p_1}(1+\mu)^{p_2}J_{x_i}(R_0)^{1/2}.
\]
It follows by a fixed number of local fourth-order Sobolev estimates and Sobolev embedding; only finitely many fixed graph derivatives are needed. Since the central boundary patch has fixed area, \eqref{m:clamped-small} gives
\[
T\le C A^{p_3-c_3N}e^{C\mu}M.
\]
Combining this with \eqref{m:clamped-anchor} and \eqref{m:clamped-cauchy-eq} gives
\[
e^{-N/2}\le C A^{\theta p_3-\theta c_3N}e^{C\mu}.
\]
Choose $A=A_C$ large enough that $\theta c_3\log A>1$, then choose $M_C$ large enough. The displayed inequality forces $N\le M_C(1+\mu)$, proving the assertion. The all-centers and all-scales quantifiers of \eqref{m:HC} are essential because witness centers can lie above the boundary and their radii can be arbitrarily small.
\end{proof}

\begin{lemma}[A bounded-frequency clamped terminal cap]\label{m:clamped-terminal}
For each fixed $N_*$, a smooth clamped cap with $\mu(Q)\le1$ and $N_C(Q)\le N_*$ satisfies
\begin{equation}\label{m:clamped-terminal-eq}
\HH^{n-1}(\Z(u)\cap W(Q))\le C_{N_*}\ell(Q)^{n-1}.
\end{equation}
The estimate holds for nonzero clamped solutions on the reserved outer cap. The constant may depend on $N_*$ and on fixed finite-order bounds for the smooth charts.
\end{lemma}
\begin{proof}
\emph{Step 1. One controlled child.}
Normalize $\ell(Q)=1$ and $\|u\|_{L^2(W_{\rm out})}=1$. Use the interior ball of Lemma~\ref{m:clamped-cauchy}. Since $W_{\rm out}\subset B_{L\rho}(y_0)$, the index bound gives
\[
 \int_{B_\rho(y_0)}|u|^2\ge e^{-N_*}J_{y_0}(L\rho)\ge e^{-N_*}.
\]
Fix a smaller connected cap containing this closed ball and the central graph patch, compactly inside the outer cap relative to the physical face. For instance take $q=\lfloor n/2\rfloor+5$ and $0<\beta<\min\{1,q-3-n/2\}$. Use uniform $C^{q+4}$ bounds for the graph and its coordinates, available from the fixed smooth atlas after rescaling. The clamped estimate \cite{ADN} on nested caps gives $\|u\circ F\|_{H^q}\le C_q$ for $\mu\le1$. Rellich and Sobolev embedding give precompactness in $C^3$ and boundedness in $C^{3,\beta}$ on a slightly smaller closed cap. The graphs have a subsequence converging at least in $C^4$. Change of variables, including the uniformly bounded positive Jacobians, passes the displayed interior mass to the limit. Thus the limit solves the limiting plate equation in physical coordinates and is nonzero.

The clamped traces imply $U_j(y,0)=\partial_tU_j(y,0)=0$. Write
\[
U_j(y,t)=t^2V_j(y,t),\qquad
V_j(y,t)=\int_0^1(1-a)\partial_t^2U_j(y,at)\,\mathrm da.
\]
The quotients are compact in $C^1$. For a limit, if $\partial_t^2U_\infty$ is nonzero at some point of a central boundary patch, $V_j$ has a fixed sign on a smaller cap for all sufficiently close data. Otherwise that second derivative vanishes on the patch. Then $\partial_t^3U_\infty$ must be nonzero somewhere on it: if it also vanished, tangential differentiation of the zero traces would give zero physical traces $u,\partial_\nu u,\Delta u,\partial_\nu\Delta u$. Zero extension would solve $(\Delta^2-\mu_\infty^4)\widetilde u=0$ across the patch, by Green's formula, and analytic unique continuation would contradict the nonzero limit.

In the second case, $\partial_tV_\infty=\partial_t^3U_\infty/6$ is nonzero at a boundary point. On a smaller rectangular cap, nearby quotients have $|\partial_tV_j|\ge c>0$ and bounded tangential gradient. Each vertical line has at most one zero, and for two zeros $(y,t)$ and $(z,b)$ the mean-value theorem gives
\[
c|t-b|\le C|y-z|.
\]
Thus the zero set lies on a uniformly Lipschitz graph over a subset of the tangential base and has bounded $(n-1)$-measure there. In both cases compactness of the normalized family supplies a uniform patch radius and measure constant: take a finite cover by neighborhoods of limiting data on which the indicated nonvanishing derivative persists. A fixed fine subdivision has one whole child in such a patch.
\emph{Step 2. The whole cap.} Repeat this good-child construction at each surviving boundary node, stopping on its good child. The local frequency remains bounded and the nested scalar index remains at most $N_*$. The band above the boundary row has bounded nodal measure by Lemma~\ref{lem:analytic}: scalar anchors are compared with the outer norm by a fixed number of admissible $D_L$ comparisons. If the good-child subdivision has $A_0^{n-1}$ children, the total remaining boundary area weight decreases by the factor $1-A_0^{-(n-1)}$. Sum this geometric series, charging each band and stopped child. Every positive-height point is charged at a finite depth. This proves \eqref{m:clamped-terminal-eq}, including artificial grid faces, without assuming a global nodal-measure bound in advance.
\end{proof}

\section{Summation, the global bound and sharpness}\label{m:sharp-section}

\begin{lemma}[Summation over boundary cubes]\label{lem:budget}
Let $d=n-1$. Above a base of side $\ell$ subdivide the base into $m=A^d$ equal children. Suppose $N\ge0$ is hereditary, $\mu\ge0$ becomes $\mu/A$ on every child, and some child has index at most $N/2$ whenever $N>M(1+\mu)$, where $M\ge1$. Suppose also that the cap is charged in full when $\mu\le1$ and $N\le2M$. If a terminal cap or the band above a nonterminal boundary row costs at most $C(1+N+\mu)\ell^d$, then the sum of all charges below a root is at most $C'(1+N_{\rm root}+\mu_{\rm root})\ell_{\rm root}^d$.
\end{lemma}
\begin{proof}
Fix $B\ge4M+4$, write $F=N+B\mu+1$ and give each node weight $F\ell^d$. At a high-index node the difference between $F$ and the mean child value is at least $N/(2m)+B(1-A^{-1})\mu$. Since $N>1$, this is at least $\min\{1/(4m),1-A^{-1}\}F$. At every other nonterminal node $\mu>1$ and $F\le(2M+B+1)\mu$, so the difference is at least $B(1-A^{-1})F/(2M+B+1)$. With $\epsilon$ the smaller contraction in these two cases,
\begin{equation}\label{eq:budget}
 \sum_qF(q)\ell(q)^d\le(1-\epsilon)F\ell^d .
\end{equation}
Terminal nodes carry no weight to the next generation. Thus generation $j$ costs at most $C(1-\epsilon)^jF_{\rm root}\ell_{\rm root}^d$. Choose closed lateral faces for both the bands and the children. Each positive-height point leaves the surviving boundary rows in finitely many steps. Countable subadditivity gives the estimate even if a nodal hypersurface lies in a grid face.
\end{proof}

\subsection{The two collars and the interior}\label{m:collars}

\begin{proposition}[The sharp Navier collar bound]\label{m:navier}
For a sufficiently small fixed $d_N>0$,
\begin{equation}\label{m:navier-total}
\HH^{n-1}\bigl(\Z(u)\cap\{\dist(x,\GN)<d_N\}\bigr)\le CK.
\end{equation}
\end{proposition}
\begin{proof}
Choose finitely many fixed buffered Navier root cubes whose physical caps cover a collar of width $d_N$. Their enlargements avoid $\GC$. By Lemma~\ref{m:navier-growth}, each root has maximal lifted index $N(Q)\le CK$. Its local frequency is $k\ell(Q)$, which need not be bounded. The arbitrary-frequency theorem Theorem~\ref{thm:local-all}, together with the exact scaling in Remark~\ref{rem:scaling}, gives
\[
\HH^{n-1}\bigl(\Z(u)\cap\Phi(Q^+)\bigr)
\le C\ell(Q)^{n-1}\bigl[1+N(Q)+k\ell(Q)\bigr]\le C_QK.
\]
Summing over this fixed finite cover proves \eqref{m:navier-total}. The entire collar is counted by the boundary recursion, so no spectral number of separate distance layers is introduced.
\end{proof}

\begin{proof}[Proof of Theorem~\ref{m:sharp}]
For each clamped node set
\[
\mathcal F_C(Q)=N_C(Q)+B_C\mu(Q)+1,
\qquad B_C\ge4M_C+4.
\]
Use Lemma~\ref{m:clamped-terminal} for nodes with $\mu\le1$ and $N_C\le2M_C$. For the other nodes, Lemma~\ref{m:clamped-high} supplies the child required in Lemma~\ref{lem:budget}, and \eqref{m:clamped-heredity} gives $N_C(q)\le N_C(Q)$ and $\mu(q)=\mu(Q)/A_C$.

To estimate a clamped band, cover its closed positive-height portion by a fixed number of balls of radius $s\asymp\ell(Q)/A_C$, with the balls $B_{L_1s}$ inside the domain. For every simplex anchor of radius $\epsilon_0s$ use its center $z\in E(Q)$ and choose $L$ initially so that $B_{L\epsilon_0s}(z)$ contains the corresponding outer analytic ball. Then the anchor has $L^2$ norm at least $e^{-N_C(Q)/2}$ times that outer norm. Formula \eqref{eq:analytic-measure} gives a total band cost at most $C\mathcal F_C(Q)\ell(Q)^{n-1}$, with $C$ allowed to depend on the fixed $A_C$. Hence the convergent boundary-tree sum gives
\[
\HH^{n-1}(\Z(u)\cap W(Q_0))
\le C\ell(Q_0)^{n-1}\bigl[N_C(Q_0)+k\ell(Q_0)+1\bigr].
\]
By Lemma~\ref{m:HC-seed}, finitely many fixed clamped roots contribute at most $CK$. The Navier collar contributes at most $CK$ by Proposition~\ref{m:navier}; its proof uses only the fixed seeds and remains valid for smooth clamped components. On the compact interior remainder, fixed-radius seeds and Lemma~\ref{lem:analytic} give the same bound. Thus the full nodal measure is at most $C_\Omega(1+k)$. The energy identity and the Dirichlet Poincar\'e inequality give $k\ge k_1(\Omega)>0$, so the additive constant is absorbed. This proves \eqref{m:sharp-bound}.
\end{proof}

\begin{remark}[Precise logical dependence]\label{m:HC-status}
A small-radius witness in Lemma~\ref{m:clamped-high} must force growth at the same center on a sequence of larger radii. The implication is \eqref{m:HC}. The global bound \eqref{m:global-L2-doubling} by itself gives only an upper bound at each radius. The two inputs from the fully clamped problem \cite{LTY} are the annular estimate \eqref{m:annular-input-eq} and the flat estimate \eqref{m:flat-input}. Theorem~\ref{m:local-doubling} combines them on a single clamped patch. Lemma~\ref{m:seeds} supplies the root masses, and Lemmas~\ref{m:clamped-high} and~\ref{m:clamped-terminal} give the two clamped subdivision steps. For the Navier collars, Proposition~\ref{m:navier} instead uses Theorem~\ref{thm:local-all} for the reflected lift.
\end{remark}

\subsection{A sharp sequence on a mixed annulus}

\begin{proposition}[Sharpness on a smooth mixed annulus]\label{m:mixed-sharpness}
On a fixed annulus $\Omega=\{a<|x|<b\}$, $0<a<b$, impose Navier conditions on the inner sphere and clamped conditions on the outer sphere. There are real eigenpairs $(\lambda_m,u_m)$ with
\[
c m^4\le\lambda_m\le C m^4,
\qquad \HH^{n-1}(\Z(u_m))\ge c m\ge c'\lambda_m^{1/4}.
\]
Thus the exponent in Theorem~\ref{m:sharp} cannot be decreased even when both boundary types are nonempty.
\end{proposition}
\begin{proof}
Choose the real spherical harmonic of degree $m$ obtained by restricting $\operatorname{Re}(x_1+i x_2)^m$ to the unit sphere, and denote it by $Y_m$. Minimize the plate Rayleigh quotient on its radial sector $R(r)Y_m(\omega)$ in the variational space \eqref{m:weak}. Compactness gives a minimizer. The angular $L^2$ projection $P_m$ onto the span of $Y_m$ at each radius preserves the form domain and commutes with $\Delta$, since $\Delta_{S^{n-1}}Y_m=-m(m+n-2)Y_m$ and the Dirichlet trace and radial differentiation commute with the angular projection. Thus $\int\Delta u\,\Delta(I-P_m)\varphi=\int u(I-P_m)\varphi=0$ for every test $\varphi\in\mathcal V$. Testing the minimizing equation on $P_m\varphi$ proves the eigenvalue equation on the full form domain. The essential radial traces are $R(a)=R(b)=R'(b)=0$; the weak equation supplies the Navier trace at $a$.

For a function in this sector, the angular Dirichlet energy gives
$\norm[L^2]{\nabla u}^2\ge m(m+n-2)b^{-2}\norm[L^2]{u}^2$.
Since $u=0$ on both boundary components, integration by parts and Cauchy--Schwarz imply
$\norm[L^2]{\Delta u}\norm[L^2]{u}\ge\norm[L^2]{\nabla u}^2$.
These inequalities prove the lower eigenvalue bound. A fixed nonzero radial test function in $C_c^\infty(a,b)$ gives the upper bound, because
\[
\Delta(RY_m)=\left(R''+\frac{n-1}{r}R'-\frac{m(m+n-2)}{r^2}R\right)Y_m.
\]
Finally the zero set of $Y_m$ contains $m$ distinct hyperplanes through the $(n-2)$-dimensional axis $x_1=x_2=0$, restricted to the sphere. Hence $\Z(u_m)$ contains their sections with the annulus. Each section has the same positive $(n-1)$-measure, and distinct sections intersect in sets of zero $(n-1)$-measure. Their sum is at least $cm$. Combining with $\lambda_m\le Cm^4$ proves the result.
\end{proof}

\par\addvspace{2\baselineskip}
\appendix
\pdfbookmark[0]{Appendices}{appendices-overview}
\begin{center}
  {\LARGE\bfseries Appendices\par}
\end{center}\nopagebreak[4]

\section{Local geometry and second-order tools}\label{sec:tools}
\noindent\textbf{Relation to the main text.}\enspace Section~\ref{m:regularity} uses the odd reflection in Lemma~\ref{lem:reflection} to obtain the interior mass bound in Lemma~\ref{m:seeds}. In Section~\ref{m:global}, the same reflection and the frequency estimate in Lemma~\ref{lem:frequency} give the doubling bound for the Navier lift in Lemma~\ref{m:navier-growth}.

All local coordinate neighborhoods have fixed buffers. Enlargements by fixed factors are made inside these buffers. Constants may depend on these factors, on $n$, and on a fixed $C^{2,\alpha}$ boundary atlas. They do not depend on a shrinking scale or on $\kappa\in[0,1]$.

\begin{lemma}[Adapted flattening and reflection]\label{lem:reflection}
A $C^{2,\alpha}$ hypersurface admits local $C^{2,\alpha}$ flattening maps $\Phi$ with uniformly bounded inverses such that the tangential columns of $D\Phi$ are orthogonal to its normal column on $\{y_n=0\}$. Set
\[
b=|\det D\Phi|,\qquad A=bD\Phi^{-1}D\Phi^{-T}.
\]
Then $A$ is symmetric uniformly elliptic, $A,b$ are one-sided $C^{1,\alpha}$, and $A_{in}=0$ on the flat boundary for $i<n$. If $f$ is harmonic and has zero Dirichlet trace on the physical boundary, the odd extension of $f\circ\Phi$ solves a full-ball divergence-form equation with Lipschitz principal coefficients. The same conclusion holds in a product with the variable $t$; the coefficient matrix before reflection is $\operatorname{diag}(A,b)$.
\end{lemma}
\begin{proof}
Write the hypersurface as $X(y')=(y',g(y'))$, with $g\in C^{2,\alpha}$, and let $\nu(y')$ be a chosen unit normal. A regularized extension of this $C^{1,\alpha}$ normal field gives
\[
\Phi(y',s)=X(y')+s\nu_s(y'),\qquad s\geq0,
\]
where $\nu_s$ is convolution with a smooth symmetric mollifier at scale $s$. At $s=0$, $D_{y'}\Phi=D_{y'}X$ and $\partial_s\Phi=\nu$. The cancellation of the first moments of the mollifier, together with the $C^{1,\alpha}$ norm of $\nu$, shows that the second derivatives of $\Phi$ extend with exponent $\alpha$. For example, $sD_{y'}^2\nu_s=O(s^\alpha)$, $D_{y'}\nu_s\to D_{y'}\nu$, and $\partial_s\nu_s+s\partial_s^2\nu_s=O(s^\alpha)$. The standard estimates for mollified H\"older functions also give the corresponding two-point H\"older bounds. Shrinking the neighborhood makes $\Phi$ a diffeomorphism. This construction avoids assuming a third derivative of $g$, as an unregularized normal-coordinate formula would do away from the boundary.

The change-of-variables formula gives $\diver(A\nabla(f\circ\Phi))=0$. Let $J=\operatorname{diag}(1,\ldots,1,-1)$ in the $y$ variables. For $y_n<0$, extend
\[
A^e(y)=JA(y',-y_n)J,\qquad b^e(y)=b(y',-y_n).
\]
The mixed entries are reflected oddly and vanish at the interface, so all entries are Lipschitz. Reflect $f\circ\Phi$ oddly. Testing with the odd part of a full-ball test function proves the reflected weak equation. The normal flux matches on the two sides; hence no interface distribution appears. Since $\Phi$ is independent of $t$, the product coefficient is block diagonal as asserted. Further spatial dilations improve the Lipschitz seminorms, and the ellipticity bounds remain uniform.
\end{proof}

For a nonzero solution $F$ of a full-ball scalar equation $\diver(B\nabla F)=0$, with $B$ symmetric, uniformly elliptic and Lipschitz, put
\[
m_z(s)=\norm[L^2(B_s(z))]{F},\qquad
\mathcal D_F(z,s)=\log_2\frac{m_z(2s)}{m_z(s)}.
\]
All such norms are positive by unique continuation.

\begin{lemma}[Frequency consequences with buffers]\label{lem:frequency}
Let $F\not\equiv0$ solve the full-ball equation with coefficients from the fixed elliptic Lipschitz class. Choose $b\ge8$ allowing for ellipticity normalization of all balls. One may use the following nested sets. Fix $z_0$ and $a>0$ so that $B_{32ba}(z_0)$ lies in the reflected coordinate patch and its radius is below the geometric frequency radius. Take $z\in B_a(z_0)$, $s\le a$. In (iii) require also $R\le a$ and $B_{bR}(z)\subset B_{32ba}(z_0)$. With constants depending only on the coefficient class and fixed ratios, the following hold.
\begingroup
\renewcommand{\labelenumi}{(\roman{enumi})}
\begin{enumerate}
\setlength{\itemsep}{4pt}
\item Fixed-ratio three-ball inequalities hold with constants and exponents determined only by the coefficient class and the ratios.
\item On a fixed buffered root region, the supremum of $\mathcal D_F(z,s)$ over centers in its smaller part and admissible smaller radii is bounded by
\[
C\left(1+\log\frac{\norm[L^2(B_{R_*})]{F}}{\norm[L^2(B_{r_*})]{F}}\right),
\]
where one may take $r_*=a/4$, $R_*=16ba$ and center both balls at $z_0$. Another fixed interior observation ball can be connected to this one by a finite overlapping chain inside the larger patch.
\item There are constants $b_0,C,c>0$ such that if $\mathcal D_F(z,s)=D$, $\rho\geq b_0s$, and $R\geq b_0\rho$, then
\begin{equation}\label{eq:lower-growth}
m_z(\rho)\leq e^{C(D+1)}(\rho/R)^{cD-C}m_z(R).
\end{equation}
The factor $b$ reserves the enlarged balls required by the ellipsoids used in normalizing $B(z)$. In applications involving the plate trace the centers are $z=(y,0)$ in spatial and lift coordinates. Frequency estimates for the full reflected equation also hold for the other centers in the displayed neighborhood.
\end{enumerate}
\endgroup
\end{lemma}
\begin{proof}
These are the buffered forms of the frequency estimates in \cite{GL86,GL87}. We recall the passage to the forms used here. After normalizing $B(z)$, the adapted frequency $\mathfrak f_z$ is nonnegative, $e^{Cs}\mathfrak f_z(s)$ is nondecreasing, and the logarithmic derivative of the associated surface height is $(d-1+2\mathfrak f_z(s))/s+O(1)$, where $d$ is the dimension of the equation. Integration compares fixed-ratio norm growth with the frequency at slightly smaller and larger radii. Integrating once more gives the three-ball inequality and the growth between separated scales. In particular, growth $D$ between $s$ and $2s$ forces $\mathfrak f_z$ at a fixed multiple of $s$ to be at least $cD-C$. Its almost monotonicity on the interval up to $R$ gives \eqref{eq:lower-growth}; fixed endpoint enlargements account for the factor $e^{C(D+1)}$.

For (ii), first compare a ball about $z$ with fixed balls about the root center. A finite chain of fixed-ratio three-ball inequalities controls the logarithmic loss in the denominator by the displayed root norm ratio. The integrated frequency formula then bounds every smaller-scale doubling index. If a radius is comparable to the root scale, the same conclusion follows directly from these fixed-ball comparisons. The number of balls is fixed geometrically. Bounded changes of the surface-height density and of the normalized metric affect only the constants and the buffers. Thus all statements are uniform for the indicated class; no parameter from the plate equation occurs in them.
\end{proof}

\begin{lemma}[Quantitative Cauchy propagation]\label{lem:cauchy}
Fix a small boundary product patch in the coordinates of Lemma~\ref{lem:reflection}. There are a smaller interior ball $B_*$, a boundary disk $\Sigma$, and $C>0$, $\theta\in(0,1)$ such that a harmonic function $H$ vanishing on the physical product boundary satisfies
\begin{equation}\label{eq:cauchy}
\norm[L^2(B_*)]{H}\leq C M^{1-\theta}\epsilon^\theta,
\quad
M=\norm[L^2(\text{buffered patch})]{H},\quad
\epsilon=\norm[L^2(\Sigma)]{\partial_\nu H}.
\end{equation}
Coordinate and physical norms, and normal and conormal data, are comparable with fixed constants.
\end{lemma}
\begin{proof}
Flatten the boundary, and extend $H$ by zero to the lower half-ball. The extended function belongs to $H^1$ and satisfies a divergence-form equation whose right-hand side is the conormal datum times surface measure on the interface. The trace theorem bounds this right-hand side in $H^{-1}$ by $C\epsilon$. The zero-Dirichlet solution $W$ of the same inhomogeneous equation on the full ball satisfies $\norm[H^1]{W}\leq C\epsilon$ by coercivity. Thus the difference between the zero extension and $W$ solves the homogeneous full-ball equation. On a fixed lower ball its norm is at most $C\epsilon$, while on the full ball it is at most $C(M+\epsilon)$. A finite chain of the three-ball inequalities in Lemma~\ref{lem:frequency} gives the bound $C(M+\epsilon)^{1-\theta}\epsilon^\theta$ on $B_*$. Add $W$. If $\epsilon\leq M$, this is \eqref{eq:cauchy}; if $\epsilon>M$, the assertion follows from the definition of $M$.
\end{proof}

\section{Exact harmonic lifting and norm transfer}\label{sec:lift}
\noindent\textbf{Relation to the main text.}\enspace The norm comparisons for the harmonic lift introduced in Section~\ref{m:regularity} are used in Lemma~\ref{m:seeds} and again in Section~\ref{m:global} to pass from lifted to scalar doubling in Corollary~\ref{m:global-scalar}. Section~\ref{m:collars} uses the scaling in Remark~\ref{rem:scaling} when applying the local Navier estimate in Proposition~\ref{m:navier}.

Let $D$ be a physical boundary patch with the geometry fixed above. We work with
\begin{equation}\label{eq:local-system}
\Delta U=V,\qquad\Delta V=\kappa^4U\quad\hbox{in }D,
\qquad U=V=0\quad\hbox{on its boundary patch}.
\end{equation}
Unless specified otherwise, $0\leq\kappa\leq1$. Companion variables are always spatially normalized: under a dilation by $r$, the pair becomes $(U(x_0+r\cdot),r^2V(x_0+r\cdot))$ and the parameter becomes $r\kappa$.

\begin{lemma}[Local energy and regularity]\label{lem:energy}
With fixed buffered enlargements, solutions of \eqref{eq:local-system} satisfy
\begin{equation}\label{eq:caccioppoli}
\norm[L^2(D\cap B_s)]{V}
\leq C(s^{-2}+\kappa^2)\norm[L^2(D\cap B_{2s})]{U}.
\end{equation}
The estimate also holds on interior balls. If $\kappa s$ is bounded, boundary regularity gives, with a fixed enlargement $B_{c_1s}$,
\begin{equation}\label{eq:gradient}
\sup_{D\cap B_s}|\nabla U|\leq Cs^{-1-n/2}\norm[L^2(D\cap B_{c_1s})]{U},
\quad
\sup_{D\cap B_s}|\nabla V|\leq Cs^{-3-n/2}\norm[L^2(D\cap B_{c_1s})]{U}.
\end{equation}
At fixed normalized scales, $U,V$ have uniform one-sided $C^{2,\alpha}$ bounds in terms of a buffered $L^2$ norm of $U$.
\end{lemma}
\begin{proof}
Use a cutoff $\eta$ with first and second derivatives bounded by $Cs^{-1}$ and $Cs^{-2}$. Set $M=\norm[L^2(D\cap B_{2s})]{U}$,
\[
A_0=\int\eta^4V^2,\qquad B_0=\int\eta^2|\nabla U|^2.
\]
Testing the first equation gives $B_0\leq CM A_0^{1/2}+Cs^{-2}M^2$. Testing the second equation against $\eta^4U$, and integrating twice, gives
\[
\kappa^4\int\eta^4U^2
=\int\eta^4V^2+8\int\eta^3V\nabla\eta\cdot\nabla U
+\int UV\Delta(\eta^4).
\]
For variational data Lemma~\ref{lem:transposition} gives the zero trace of $V$ and $H^2$ regularity before either integration by parts. Both physical boundary terms therefore vanish. Young's inequality yields $A_0\leq C(s^{-4}+\kappa^4)M^2$. Nested cutoffs give \eqref{eq:caccioppoli} with the stated radii. Boundary and interior Dirichlet estimates, applied successively to the two second-order equations, give $W^{2,p}$ bounds and then $C^{2,\alpha}$ bounds. In arbitrary dimension one iterates the Sobolev step finitely many times before applying Schauder estimates. Scaling gives \eqref{eq:gradient}. These are ordinary second-order estimates; see \cite{ADN}.
\end{proof}

Define
\begin{equation}\label{eq:CS}
C_\kappa(t)=\frac{\cosh(\kappa t)+\cos(\kappa t)}2,
\qquad
S_\kappa(t)=\begin{cases}
\dfrac{\cosh(\kappa t)-\cos(\kappa t)}{2\kappa^2},&\kappa>0,\\[4pt]
t^2/2,&\kappa=0.
\end{cases}
\end{equation}

\begin{lemma}[Harmonic lift and its stable trace]\label{lem:lift}
The real function
\begin{equation}\label{eq:lift}
H(x,t)=C_\kappa(t)U(x)-S_\kappa(t)V(x)
\end{equation}
is harmonic in the product patch and vanishes on the Navier product boundary. Also $H(x,0)=U(x)$ and $H_{tt}(x,0)=-V(x)$. For every fixed $M_*>0$ there are positive constants $c,C$ such that, whenever $\kappa s\leq M_*$,
\begin{equation}\label{eq:gram}
c s(|U|^2+s^4|V|^2)
\leq\int_{-s}^s|H(x,t)|^2\dd t
\leq C s(|U|^2+s^4|V|^2).
\end{equation}
All trace comparisons use a spatial center $x_0$ and the lifted center $(x_0,0)$, with each indicated larger ball contained in the reserved chart (or in its reflection). In particular, for fixed buffer factors and $\kappa s\leq M_*$,
\begin{equation}\label{eq:trace}
c s^{1/2}\norm[L^2(D\cap B_{c_2s})]{U}
\leq\norm[L^2((D\times\R)\cap B_s)]{H}
\leq C s^{1/2}\norm[L^2(D\cap B_{c_3s})]{U}.
\end{equation}
Here $c_2>0$ is small and $c_3>1$ is fixed. These comparisons persist in adapted reflected coordinates.
\end{lemma}
\begin{proof}
The identities $C_\kappa''=\kappa^4S_\kappa$ and $S_\kappa''=C_\kappa$ prove harmonicity directly using \eqref{eq:local-system}. At $\kappa=0$, the lift is $U-t^2V/2$, so there is no singular normalization.

Put $t=s\tau$ and $a=\kappa s$. Since $C_\kappa(s\tau)=C_a(\tau)$ and $S_\kappa(s\tau)=s^2S_a(\tau)$, the integral in \eqref{eq:gram} is the quadratic form of the Gram matrix of $C_a,-S_a$ on $(-1,1)$, applied to $(U,s^2V)$, multiplied by $s$. This matrix depends continuously on $a\in[0,M_*]$ and is positive definite at every $a$: evaluation at zero and taking the second derivative at zero show that the two functions are linearly independent. Compactness proves the uniform upper and lower eigenvalue bounds, including at $a=0$.
Explicitly, the limiting matrix is
\[
G_0=\begin{pmatrix}2&-1/3\\-1/3&1/10\end{pmatrix},\qquad
\det G_0=4/45>0.
\]

A ball in $n+1$ variables is contained in a product cylinder and contains a smaller product cylinder. Integrating \eqref{eq:gram} and using \eqref{eq:caccioppoli} proves \eqref{eq:trace}. Coordinate Jacobians and ball comparisons are uniformly bounded. Odd reflection preserves the two-term formula in $t$, and replaces half-ball norms by comparable full-ball norms.
\end{proof}

\begin{remark}\label{rem:scaling}
The lift is compatible with spatial rescaling. If $U_r(y)=U(x_0+ry)$, $V_r(y)=r^2V(x_0+ry)$ and $\kappa_r=r\kappa$, their lift is exactly $H(x_0+ry,rt)$. Thus a single lifted function can drive every generation of a boundary recursion. This compatibility is essential; independently renormalizing a spectral companion at each generation would not give the same maximal index.
\end{remark}

In both norm transfers below the lifted balls are again centered at $(x_0,0)$ and the spatial balls at $x_0$. For a center on the reflected side interpret the spatial comparison by reflection back to the physical side. The energy estimate \eqref{eq:caccioppoli} holds for every $\kappa\ge0$, as its proof shows. We now allow $\kappa$ to be arbitrary and need an estimate when $\kappa s$ is not bounded. The elementary bounds $|C_\kappa(t)|\leq Ce^{\kappa|t|}$ and $|S_\kappa(t)|\leq Ct^2e^{\kappa|t|}$, together with \eqref{eq:caccioppoli}, give
\begin{equation}\label{eq:large-parameter-upper}
\norm[L^2(B_s^{n+1})]{H}
\leq Cs^{1/2}e^{C\kappa s}(1+(\kappa s)^2)
\norm[L^2(D\cap B_{c_3s})]{U}.
\end{equation}
Conversely, the trace inequality and the Caccioppoli inequality for the scalar harmonic equation imply
\begin{equation}\label{eq:large-parameter-lower}
\norm[L^2(D\cap B_s)]{U}
\leq Cs^{-1/2}\norm[L^2(B_{c_4s}^{n+1})]{H}.
\end{equation}
For a boundary patch this is proved after reflection, using the uniformly elliptic Lipschitz equation for the lift. Indeed, the squared trace norm is bounded by $C(s^{-1}\norm{H}^2+s\norm{\nabla H}^2)$ on a surrounding cylinder, and elliptic Caccioppoli controls the gradient. There is no spectral coefficient in this lifted equation.

\begin{lemma}[Two boundary traces from the lift at arbitrary frequency]\label{lem:lift-traces}
Let $\kappa\ge0$ in \eqref{eq:local-system}, and let $H$ be its harmonic lift. On nested product neighborhoods of a Navier boundary piece $\Gamma_\rho$, of diameter comparable to $\rho$, one has
\begin{equation}\label{eq:lift-traces}
\sup_{\Gamma_\rho}\bigl(|\partial_\nu U|+\rho^2|\partial_\nu V|\bigr)
\le C\rho^{-(n+3)/2}\norm[L^2(\mathcal B_{c\rho}^+)]{H}.
\end{equation}
Here $\mathcal B_{c\rho}^+$ is a fixed buffered product cap, and the constant is independent of $\kappa$. The right side can equivalently be replaced, up to fixed constants and buffers, by the norm of the reflected lift on a full coordinate ball.
\end{lemma}
\begin{proof}
The physical product boundary is independent of $t$. Difference quotients in $t$ preserve both the harmonic equation and its homogeneous Dirichlet condition. Iterating the local energy inequality on nested product caps therefore gives
\[
\norm[L^2(\mathcal B_{c_1\rho}^+)]{\partial_t^2H}
\le C\rho^{-2}\norm[L^2(\mathcal B_{c\rho}^+)]{H}.
\]
This assertion uses only the harmonic equation, not the size of $\kappa$ in the explicit formula for $H$. For a harmonic function $W$ with zero trace on a $C^{2,\alpha}$ boundary, the scaled local boundary gradient estimate is
\[
\sup_{\mathcal B_{c_0\rho}^+}|\nabla W|
\le C\rho^{-1-(n+1)/2}\norm[L^2(\mathcal B_{c_1\rho}^+)]{W}.
\]
Apply it first to $H$ and then to $\partial_t^2H$. The latter has zero boundary trace as well. Since $H(x,0)=U(x)$ and $H_{tt}(x,0)=-V(x)$, their spatial normal derivatives at $t=0$ are precisely the two traces in \eqref{eq:lift-traces}. The adapted coordinate maps and their Jacobians are uniformly bounded, proving the reflected version. All buffers are fixed before any subdivision.
\end{proof}

\section{A controlled boundary child at bounded index}\label{sec:low}
\noindent\textbf{Relation to the main text.}\enspace This appendix supplies the bounded-index step behind the Navier collar estimate in Section~\ref{m:collars}. In Theorem~\ref{thm:local}, Proposition~\ref{prop:low} bounds the nodal measure in one boundary child. That bounded-frequency theorem controls the terminal cubes in Theorem~\ref{thm:local-all}, which Section~\ref{m:collars} applies in the proof of Proposition~\ref{m:navier}.

The following argument is in physical coordinates and uses only $C^{2,\alpha}$ boundary regularity. Let $h$ be positive harmonic in a buffered boundary patch, zero on its boundary portion, and normalized at a fixed interior point. One can construct it by prescribing nonnegative data on the artificial boundary, separated from the portion under study. Boundary estimates and the Hopf lemma give uniform upper bounds and
\[
c\dist(x,\partial D)\leq h(x)\leq C\dist(x,\partial D),\qquad
|\partial_\nu h|\geq c
\]
on smaller patches. The chosen normal orientation is immaterial below.

\begin{lemma}[Regularity of the quotients]\label{lem:quotient}
Suppose the buffered norm of a solution of \eqref{eq:local-system} is bounded and $0\leq\kappa\leq1$. Then, on a smaller closed boundary patch,
\begin{equation}\label{eq:quotient-bound}
q=U/h,\quad p=V/h,\qquad
\norm[C^{2,\alpha}]{q}+\norm[C^{2,\alpha}]{p}\leq C.
\end{equation}
Furthermore,
\begin{equation}\label{eq:quotient-eqn}
h\Delta q+2\nabla h\cdot\nabla q=hp,
\qquad \partial_\nu q=0\quad\hbox{on the boundary}.
\end{equation}
If $q$ vanishes identically on an open boundary portion, then
\begin{equation}\label{eq:second-normal}
3\partial_{\nu\nu}q=p\quad\hbox{there}.
\end{equation}
\end{lemma}
\begin{proof}
By Lemma~\ref{lem:energy}, $U,V$ have uniform $C^{2,\alpha}$ bounds. Apply the inhomogeneous higher-order boundary Harnack estimate \cite[Theorem 3.1, with regularity order $2$]{DS} to the positive harmonic denominator $h$ and to $U$ with forcing $V$, and then to $V$ with forcing $\kappa^4U$. The required $C^{1,\alpha}$ forcing norms are bounded. This proves \eqref{eq:quotient-bound}. Alternatively, the homogeneous harmonic quotient theorem can be applied in the product domain to the lifts of $U$ and $V$, with denominator $h(x)$ independent of $t$.

Expanding $\Delta(hq)=hp$ gives the first equation in \eqref{eq:quotient-eqn}. Taking its continuous boundary limit gives $2\partial_\nu h\,\partial_\nu q=0$, proving the Neumann condition. If also $q=0$ along that portion, its whole gradient vanishes there. Differentiating these trace identities tangentially shows that its Hessian has only a normal--normal component at a boundary point $x_0$. Write $a=\partial_\nu h(x_0)\ne0$ and $b=\partial_{\nu\nu}q(x_0)$, with $\nu$ taken inward for the expansion. Along $x_0+s\nu$, as $s\downarrow0$,
\[
h=as+o(s),\quad\nabla h=a\nu+o(1),\quad
\nabla q=bs\nu+o(s),\quad\Delta q=b+o(1).
\]
The left side of the equation for $q$ is $3ab s+o(s)$ and the right side is $ap(x_0)s+o(s)$. Division by $as$ proves \eqref{eq:second-normal}. This computation requires only the established second derivatives.
\end{proof}

\begin{lemma}[Directional line counting]\label{lem:counting}
Suppose a real function $q$, analytic in a domain and $C^2$ up to a graph boundary, satisfies $|D_e^2q|\geq c>0$ in a boundary neighborhood for every $e$ in a fixed open spherical cap $E$. Assume these directions are transverse to the graph and the neighborhood is a ball intersected with its epigraph. Then
\[
\HH^{n-1}(\{q=0\}\cap B_\rho\cap D)\leq C\rho^{n-1}.
\]
\end{lemma}
\begin{proof}
Choose the cap sufficiently transverse that, along every relevant line, the signed graph height is strictly monotone. Its intersection with the ball and epigraph is then an interval. Rolle's theorem gives at most two zeros on each such interval. The underlying zero set of a nonzero real-analytic function is countably rectifiable in dimension $n-1$, with lower-dimensional strata contributing zero to $\HH^{n-1}$. For its almost-everywhere unit normal $\nu_Z$, the area formula for projection onto $e^\perp$ gives
\[
\int_{\{q=0\}\cap B_\rho\cap D}|\nu_Z\cdot e|\dd\HH^{n-1}
\leq 2\omega_{n-1}\rho^{n-1}.
\]
Integrate over $E$. Compactness of the unit sphere and openness of $E$ imply $\inf_{|\nu|=1}\int_E|\nu\cdot e|\dd\sigma(e)>0$. This proves the result, first on compact subsets of the interior and then by exhaustion. Finiteness of the measure was not assumed in taking this limit.
\end{proof}

We use boundary cubes in flattened $y$ coordinates. A root cube $Q$ has center on $y_n=0$, its upper half is denoted $Q^+$, and its physical image is $\Phi(Q^+)$. A subdivision by $A$ partitions the $(n-1)$ tangential directions into $A^{n-1}$ cells and places a boundary cube of side $\ell(Q)/A$ over each cell.

\begin{proposition}[Bounded-index boundary child]\label{prop:low}
For each fixed $N_0$, there are an integer $A_0\geq2$ and a constant $C_0$ with the following property. Suppose the maximal lifted index on a buffered boundary cube $Q$ is at most $N_0$, and $\kappa\ell(Q)\leq1$. Then one child of its $A_0$-subdivision satisfies
\begin{equation}\label{eq:low-child}
\HH^{n-1}(\Z(U)\cap\Phi(q^+))\leq C_0\ell(Q)^{n-1}.
\end{equation}
The constants are uniform over further dilations of the fixed boundary atlas.
\end{proposition}
\begin{proof}
Normalize $\ell(Q)=1$ and the scalar $L^2$ norm on a fixed inner ball to one. By Lemma~\ref{lem:lift}, a smaller lifted ball has bounded norm. A fixed number of doublings, all controlled by $N_0$, bounds the lift on every required buffer. The lower part of \eqref{eq:gram} then bounds both $U$ and $V$ on smaller physical neighborhoods. Thus Lemmas~\ref{lem:energy} and~\ref{lem:quotient} give uniform $C^{2,\alpha}$ bounds for $U,V,h,q,p$.

After pulling back by their graph charts, these normalized data form a precompact family in $C^{2,\beta}$ on smaller patches, for any $0<\beta<\alpha$. This remains true when the physical scale and boundary point vary: the graph functions have uniform $C^{2,\alpha}$ norms, the parameters belong to $[0,1]$, and the positive functions $h$ have uniform Hopf lower bounds. Every limit solves the same system on its limiting $C^{2,\beta}$ graph domain. The scalar normalization guarantees that the limiting $U$ is nonzero.

Fix such a limit. If the boundary trace of $q$ is nonzero at a point of the central boundary patch, $q$ has a fixed sign in a smaller one-sided neighborhood. The same is true for all sufficiently close data, so that neighborhood has no interior zeros.

Otherwise $q$ vanishes on the whole central boundary patch. If $p$ also vanished there, $U,V$ would have zero values and zero normal derivatives there, since $U=hq$, $V=hp$. Their harmonic lift would have zero Cauchy data on an open product boundary patch. Extending it by zero across that patch gives a distributionally harmonic function that vanishes on an open set; Weyl's lemma and harmonic unique continuation would force the lift, and then $U$, to vanish. This contradicts the normalization. Therefore $p(x_0)\ne0$ at some central boundary point. By \eqref{eq:second-normal}, $D^2q(x_0)$ has a nonzero normal--normal entry and no other entries in a normal--tangential frame. On a sufficiently small neighborhood, $D_e^2q$ is bounded away from zero for every direction in a cap around the normal. This property persists under the $C^2$ convergence of the quotients and the boundary charts. Lemma~\ref{lem:counting} provides a uniform local measure bound for all nearby data. The quotients are analytic inside the domain, since both $U$ and the positive harmonic function $h$ are analytic there.

Each limiting datum thus has an open neighborhood in the compact family on which one physical boundary patch enjoys a uniform nodal-measure bound. A finite subcover gives a positive lower bound for the available patch radii and a common measure constant. Choose $A_0$ large enough that one whole boundary child fits in each such patch. This proves \eqref{eq:low-child}; bounded Jacobians compare coordinate and physical measures, and rescaling gives $\ell(Q)^{n-1}$.
\end{proof}

\section{High-index reduction and the local measure theorem}\label{sec:recursion}
\noindent\textbf{Relation to the main text.}\enspace Section~\ref{m:collars} applies Theorem~\ref{thm:local-all} on each Navier boundary patch and sums the resulting estimates to prove Proposition~\ref{m:navier}. This supplies the Navier collar contribution to the proof of Theorem~\ref{m:sharp}.

Let $\kappa\ge0$ and let $F(y,t)$ be the odd reflected lift in an adapted chart, with $y=(y',y_n)\in\R^n$. For a closed tangential base $I$ of side $\ell$ set
\[
 Q=I\times[-\ell/2,\ell/2],\qquad Q^+=Q\cap\{y_n>0\},\qquad
 E_N(Q)=2Q=(2I)\times[-\ell,\ell].
\]
The boundary base is $I\times\{0\}$ and its physical cap is $\Phi(Q^+)$. A boundary child uses one of the $A^{n-1}$ subbases $I_i$ and is the full cube $q_i=I_i\times[-\ell/(2A),\ell/(2A)]$. Doubling $2Q$ here is concentric dilation in all $n$ spatial variables. In particular $q_i\subset Q$ and $2q_i\subset2Q$. Fix a geometric $L$ large enough for the root and interior anchor comparisons. For the center $y_Q$ of a root take $a_Q=(2L+4\sqrt n)\ell(Q)$ and require the lifted ball $B_{32ba_Q}((y_Q,0))$ from Lemma~\ref{lem:frequency} to lie inside the full reflected patch and its uniform frequency radius. This is accomplished by shrinking each root after $L$ is fixed. For every descendant, $y\in E_N(q)$ and $s\le L\ell(q)$ give $(y,0)\in B_{a_Q}((y_Q,0))$ and $2s\le a_Q$. Define
\begin{equation}\label{eq:max-index}
N(Q)=\sup\{\mathcal D_F((y,0),s):y\in2Q,\ 0<s\leq L\ell(Q)\}.
\end{equation}
Here every norm in the index is taken in $n+1$ dimensions. The witnesses may lie on either side of the reflection plane in the spatial variable. The target center in Lemma~\ref{lem:cauchy} can be chosen in the upper central cap at lift height $t=0$ by centering its product patch on that slice. For each band ball, choose the simplex anchors close enough to its center that their spatial centers lie in $2Q$. Cover the closed band by such balls, shrinking their radii by a fixed factor to keep the analytic enlargements in the physical interior. Root, target and anchor comparisons then all use centers in $E_N(Q)$ and radii at most $L\ell(Q)$.

In particular, for every boundary descendant $q$,
\begin{equation}\label{eq:heredity}
N(q)\leq N(Q).
\end{equation}
By Lemmas~\ref{lem:frequency} and~\ref{lem:lift}, the index is finite and, on a normalized root,
\begin{equation}\label{eq:root}
N(Q)\leq C(1+D_U+\kappa),\qquad
D_U=\log\frac{\norm[L^2(D\cap B_{R_*})]{U}}{\norm[L^2(D\cap B_1)]{U}},
\end{equation}
where the outer radius is a fixed sufficiently large buffer. Bound the lifted outer norm above by the scalar outer norm using \eqref{eq:large-parameter-upper}, and bound an inner lifted norm below by the scalar inner norm using \eqref{eq:large-parameter-lower}; choose the nested radii so that every required spatial buffer is included in $B_{R_*}$. Apply part (ii) of Lemma~\ref{lem:frequency}. The logarithmic conversion cost is at most $C[1+\kappa+\log(1+\kappa)]$, which is bounded by $C(1+\kappa)$. When $\kappa\le1$ this reduces to $C(1+D_U)$.

\begin{proposition}[High-index boundary child at arbitrary frequency]\label{prop:high}
There are an integer $A\ge2$ and a constant $M_1\ge1$ such that, for every $\kappa\ge0$,
\begin{equation}\label{eq:high-child}
N(Q)>M_1\bigl(1+\kappa\ell(Q)\bigr)
\quad\Longrightarrow\quad
N(q)\le\tfrac12N(Q)
\quad\text{for at least one boundary child }q.
\end{equation}
The constants are uniform in the scale, generation, and frequency.
\end{proposition}
\begin{proof}
Normalize $\ell(Q)=1$, replacing the parameter by its rescaled value, still denoted $\kappa$, and write $N=N(Q)$. Fix a small central boundary disk $I$ and a product disk $\Sigma=I\times(-\delta,\delta)$ for Lemma~\ref{lem:cauchy}. All their buffers and comparison balls are fixed before $A$. Choose the index buffer $L$ large enough for these comparisons. Let $M$ be the reflected lifted norm on a fixed central ball containing the comparison balls of radius $\delta$ above $I$.

Suppose every boundary child has index greater than $N/2$. For each child meeting $I$, choose a witness $y_i\in2q_i$, $s_i\le L/A$ with
\[
N/2<\mathcal D_F((y_i,0),s_i)\le N.
\]
Set $\rho=B/A$, where $B$ is fixed large enough to contain the buffered child neighborhoods and to satisfy the hypotheses of Lemma~\ref{lem:frequency}. Increase $A$ until $\delta\ge b_0\rho$. The growth estimate \eqref{eq:lower-growth}, whose equation contains no spectral parameter, yields
\begin{equation}\label{eq:small-lift}
\norm[L^2(B_\rho(y_i,0))]{F}
\le C e^{C_1(N+1)}A^{-c_1N+C_1}M.
\end{equation}
Fixed changes in the radius of this ball are absorbed by the initial choice of $B$.

By Lemma~\ref{lem:lift-traces}, the two spatial normal traces on the child's physical boundary piece satisfy
\begin{equation}\label{eq:small-normal}
|\partial_\nu U|+|\partial_\nu V|
\le C\rho^{-P}\norm[L^2(B_\rho(y_i,0))]{F},
\end{equation}
with a fixed finite $P$; for example $P=(n+7)/2$ suffices after fixing the buffers and taking $\rho<1$. This estimate is valid even when $\kappa\rho>1$.

On the fixed product boundary,
\[
\partial_\nu H(x,t)=C_\kappa(t)\partial_\nu U(x)-S_\kappa(t)\partial_\nu V(x).
\]
For $|t|\le\delta$ the elementary bounds
$|C_\kappa(t)|\le Ce^{\kappa\delta}$ and
$|S_\kappa(t)|\le Ct^2e^{\kappa\delta}$ hold also at $\kappa=0$ by continuity. The children cover $I$, whose area is fixed. Consequently
\begin{equation}\label{eq:small-cauchy}
\norm[L^2(\Sigma)]{\partial_\nu H}
\le C A^{p-c_1N}e^{C_2(N+1)+C_2\kappa}M.
\end{equation}
The factor $e^{C_2\kappa}$ is retained explicitly; no bounded-frequency hypothesis has been used. The two-term formula in $t$ gives small Cauchy data on a product disk, rather than only on its slice $t=0$.

A fixed number of admissible doublings bounds the norm on the larger Cauchy patch by $e^{C_3(N+1)}M$. Applying Lemma~\ref{lem:cauchy} gives on its fixed interior target ball
\begin{equation}\label{eq:cauchy-upper}
\norm[L^2(B_*)]{F}
\le C A^{p\theta-c_1\theta N}e^{C_4(N+1)+C_4\kappa}M.
\end{equation}
The same index, using a fixed enlargement about the target center, gives
\begin{equation}\label{eq:cauchy-lower}
\norm[L^2(B_*)]{F}\ge e^{-C_5(N+1)}M.
\end{equation}
Thus
\[
\bigl(c_1\theta\log A-C_4-C_5\bigr)N
\le C_6(1+\kappa)+p\theta\log A.
\]
First fix $A$ so the coefficient on the left is positive, then fix $M_1$ large enough to make this impossible for $N>M_1(1+\kappa)$. All constants have been chosen independently of the solution and frequency. This proves \eqref{eq:high-child}.
\end{proof}

\begin{theorem}[Local Navier measure estimate at bounded frequency]\label{thm:local}
For a nonzero real solution of \eqref{eq:local-system} in a connected buffered $C^{2,\alpha}$ boundary patch, with $0\leq\kappa\ell(Q)\leq1$,
\begin{equation}\label{eq:local-measure}
\HH^{n-1}(\Z(U)\cap\Phi(Q^+))
\leq C\ell(Q)^{n-1}(1+N(Q)).
\end{equation}
In particular, on a normalized root this is bounded by $C(1+D_U)$ with $D_U$ from \eqref{eq:root}.
\end{theorem}
\begin{proof}
Fix $A,M_1$ from Proposition~\ref{prop:high}, put $N_0=2M_1$, and then fix $A_0,C_0$ from Proposition~\ref{prop:low}. Since the local frequency is at most one, $N>N_0$ implies the hypothesis of Proposition~\ref{prop:high}. Give a surviving boundary cube $Q'$ the weight
\[
w(Q')=(N(Q')+N_0)\ell(Q')^{n-1}.
\]
If $N(Q')>N_0$, subdivide by $A$ and retain every boundary child. Heredity and \eqref{eq:high-child} give
\begin{align*}
\sum_qw(q)
&\leq\left(N(Q')+N_0-\frac{N(Q')}{2A^{n-1}}\right)\ell(Q')^{n-1}\\
&\leq\left(1-\frac1{4A^{n-1}}\right)w(Q').
\end{align*}
If $N(Q')\leq N_0$, subdivide by $A_0$, count the nodal measure in the good child using \eqref{eq:low-child}, and remove it. The remaining children have total weight at most $(1-A_0^{-(n-1)})w(Q')$. The removed child's cost is at most $Cw(Q')$.

At either step, the part of $(Q')^+$ above the boundary row is a positive distance, comparable to $\ell(Q')/B$, from the boundary, where $B$ is the fixed integer $A$ or $A_0$. Cover its physical image by at most $C_B$ balls of radius comparable to $\ell(Q')/B$, chosen small enough that the enlargements needed in Lemma~\ref{lem:analytic} stay inside the physical domain. For every anchor, \eqref{eq:trace} compares its scalar norm with lifted norms at fixed comparable radii. A fixed number of doublings with centers in $2Q'$ gives the anchor lower bounds with $D\leq C_B(1+N(Q'))$. All relevant spectral products are bounded because $\kappa\ell(Q')\leq1$. By Lemma~\ref{lem:analytic}, the cost of the whole interior part is at most $C_B\ell(Q')^{n-1}(1+N(Q'))\leq Cw(Q')$. The finite choices of $B$ were fixed before the recursion, so this $C$ is independent of its generation.

With $\eta=\min\{(4A^{n-1})^{-1},A_0^{-(n-1)}\}>0$, the surviving weight of generation $j$ is at most $(1-\eta)^jw(Q)$. All charged measures sum to at most $C\eta^{-1}w(Q)$. Every compact subset of $\Phi(Q^+)$ leaves the surviving boundary rows after finitely many generations, since their heights tend uniformly to zero. Exhausting the open half-cube proves \eqref{eq:local-measure} without assuming finiteness of a residual nodal measure. Finally use \eqref{eq:root}.
\end{proof}

\begin{theorem}[Local Navier measure estimate at arbitrary frequency]\label{thm:local-all}
Let $(U,V)$ be a nonzero real solution of \eqref{eq:local-system} on a connected buffered $C^{2,\alpha}$ boundary patch, with $\kappa\ge0$. Then
\begin{equation}\label{eq:local-all}
\HH^{n-1}\bigl(\Z(U)\cap\Phi(Q^+)\bigr)
\le C\ell(Q)^{n-1}\bigl[1+N(Q)+\kappa\ell(Q)\bigr].
\end{equation}
In particular, on a unit root the right side is at most $C(1+D_U+\kappa)$, with $D_U$ as in \eqref{eq:root}. All constants are independent of $\kappa$.
\end{theorem}
\begin{proof}
Write $d=n-1$, $\mu(Q')=\kappa\ell(Q')$, and fix $A,M_1$ from Proposition~\ref{prop:high}. Choose $B_1\ge4M_1+4$ and set
\[
\mathcal F(Q')=N(Q')+B_1\mu(Q')+1,
\qquad w(Q')=\mathcal F(Q')\ell(Q')^d.
\]
A node is terminal if $\mu(Q')\le1$ and $N(Q')\le2M_1$. For such a node, the bounded-frequency theorem Theorem~\ref{thm:local} charges its entire cap at most $Cw(Q')$.

The index heredity \eqref{eq:heredity}, the scaling $\mu(q)=\mu(Q')/A$ and the high-index child from Proposition~\ref{prop:high} meet the hypotheses of Lemma~\ref{lem:budget} with $M=M_1$. The contraction \eqref{eq:budget} is uniform in the generation and the local frequency. It remains to verify the cost of the band.

The band above the boundary row is covered by a fixed number of interior balls of radii comparable to $\ell(Q')/A$, with their analytic buffers inside the physical domain. We detail the frequency dependence in the anchor comparison. For any such ball and each of the simplex anchors, \eqref{eq:large-parameter-upper}--\eqref{eq:large-parameter-lower}, and a fixed number of doublings of the single reflected lift, imply
\[
\norm[L^2(\text{anchor})]{U}
\ge \exp\{-C_A[1+N(Q')+\mu(Q')]\}
\norm[L^2(\text{outer analytic ball})]{U}.
\]
Indeed the outer scalar norm is bounded by a larger lifted norm by the harmonic trace estimate. The smaller lifted norm is bounded above by the corresponding scalar anchor norm times $e^{C\kappa r}(1+(\kappa r)^2)$; all scalar-to-lift buffers are placed inside that anchor. The logarithm of this factor is at most $C_A(1+\mu(Q'))$. Fixed-radius comparisons of the lifted balls cost at most $C_A(1+N(Q'))$. Radius powers cancel between comparable balls. These observations justify the displayed anchor bound for arbitrary $\mu(Q')$.

The interior estimate Lemma~\ref{lem:analytic} now charges the entire band at most
$C_A[1+N(Q')+\mu(Q')]\ell(Q')^d\le Cw(Q')$.
By \eqref{eq:budget}, all band and terminal charges sum to at most
$C\sum_{j\ge0}(1-\epsilon)^jw(Q)$.
Every point at positive boundary height enters a band or a terminal cap after finitely many generations. Artificial grid faces are included in the closed bands and children; multiple counting is harmless. Thus countable subadditivity gives \eqref{eq:local-all}, without assuming finite measure for a residual set. The last assertion follows from \eqref{eq:root}.
\end{proof}

\end{document}